\documentclass[a4paper,11pt]{amsart}
\usepackage{mathrsfs}
\usepackage{cases}
\usepackage{epic}
\usepackage{amsfonts}
\usepackage{graphicx}
\usepackage{amsmath}
\usepackage{amssymb, upgreek}
\usepackage{bm}
\usepackage{latexsym,todonotes}
\usepackage{pdflscape}
\usepackage[all]{xypic}
\usepackage[all]{xy}
\usepackage{color}
\usepackage{colordvi}
\usepackage{multicol}
\usepackage[linktocpage=true]{hyperref}
\hypersetup{colorlinks,linkcolor=blue,urlcolor=cyan,citecolor=blue}
\newcommand{\bvs}{\mathbf{\varsigma}}
\newcommand{\vs}{\varsigma}
\allowdisplaybreaks

\begin{document}
	\input xy
	\xyoption{all}

	\newtheorem{innercustomthm}{{\bf Main~Theorem}}
	\newenvironment{customthm}[1]
	{\renewcommand\theinnercustomthm{#1}\innercustomthm}
	{\endinnercustomthm}
	
	\newtheorem{innercustomcor}{{\bf Corollary}}
	\newenvironment{customcor}[1]
	{\renewcommand\theinnercustomcor{#1}\innercustomcor}
	{\endinnercustomthm}
	
	\newtheorem{innercustomprop}{{\bf Proposition}}
	\newenvironment{customprop}[1]
	{\renewcommand\theinnercustomprop{#1}\innercustomprop}
	{\endinnercustomthm}

	\newcommand{\iadd}{\operatorname{iadd}\nolimits}
	\newcommand{\Gr}{\operatorname{Gr}\nolimits}
	\newcommand{\FGS}{\operatorname{FGS}\nolimits}
	
	\renewcommand{\mod}{\operatorname{mod^{\rm nil}}\nolimits}
	\newcommand{\proj}{\operatorname{proj}\nolimits}
	\newcommand{\inj}{\operatorname{inj.}\nolimits}
	\newcommand{\rad}{\operatorname{rad}\nolimits}
	\newcommand{\Span}{\operatorname{Span}\nolimits}
	\newcommand{\soc}{\operatorname{soc}\nolimits}
	\newcommand{\ind}{\operatorname{inj.dim}\nolimits}
	\newcommand{\Ginj}{\operatorname{Ginj}\nolimits}
	\newcommand{\res}{\operatorname{res}\nolimits}
	\newcommand{\np}{\operatorname{np}\nolimits}
	\newcommand{\Fac}{\operatorname{Fac}\nolimits}
	\newcommand{\Aut}{\operatorname{Aut}\nolimits}
	\newcommand{\DTr}{\operatorname{DTr}\nolimits}
	\newcommand{\TrD}{\operatorname{TrD}\nolimits}
	
	\newcommand{\Mod}{\operatorname{Mod}\nolimits}
	\newcommand{\R}{\operatorname{R}\nolimits}
	\newcommand{\End}{\operatorname{End}\nolimits}
	\newcommand{\lf}{\operatorname{l.f.}\nolimits}
	\newcommand{\Iso}{\operatorname{Iso}\nolimits}
	\newcommand{\aut}{\operatorname{Aut}\nolimits}
	\newcommand{\Ui}{{\mathbf U}^\imath}
	\newcommand{\UU}{{\mathbf U}\otimes {\mathbf U}}
	\newcommand{\UUi}{(\UU)^\imath}
	\newcommand{\tUU}{{\tU}\otimes {\tU}}
	\newcommand{\tUUi}{(\tUU)^\imath}
	\newcommand{\tUi}{\widetilde{{\mathbf U}}^\imath}
	\newcommand{\sqq}{{\bf v}}
	\newcommand{\sqvs}{\sqrt{\vs}}
	\newcommand{\dbl}{\operatorname{dbl}\nolimits}
	\newcommand{\swa}{\operatorname{swap}\nolimits}
	\newcommand{\Gp}{\operatorname{Gp}\nolimits}
	
	\newcommand{\U}{{\mathbf U}}
	\newcommand{\tU}{\widetilde{\mathbf U}}
	\newcommand{\fgm}{{\rm mod}^{{\rm fg}}}
	\newcommand{\fgmz}{\mathrm{mod}^{{\rm fg},\Z}}
	\newcommand{\fdmz}{\mathrm{mod}^{{\rm nil},\Z}}
	
	\newcommand{\ov}{\overline}
	\newcommand{\und}{\underline}
	\newcommand{\tk}{\widetilde{k}}
	\newcommand{\tK}{\widetilde{K}}
	\newcommand{\tTT}{\operatorname{\widetilde{\texttt{\rm T}}}\nolimits}
		\def \cI{{\mathcal I}}
			\def \cJ{{\mathcal J}}
			\def \bbK{\mathbb{K}}
	\def\bfk{\mathbf{k}}
	\newcommand{\tH}{\operatorname{{\ch}_{\rm{tw}}}\nolimits}
	
	\newcommand{\utM}{\operatorname{\cm\ch}\nolimits}
	\newcommand{\tM}{\operatorname{\cs\cd\widetilde{\ch}}\nolimits}
	\newcommand{\rM}{\operatorname{\cm\ch_{\rm{red}}}\nolimits}
	\newcommand{\utMH}{\cs\cd\ch(\Lambda^\imath)}
	\newcommand{\tMH}{\cs\cd\widetilde{\ch}(\Lambda^\imath)}
	\newcommand{\tCMH}{{\cc\widetilde{\ch}(\K Q,\btau)}}
	
	\newcommand{\rMH}{\operatorname{\cs\cd\ch_{\rm{red}}(\Lambda^\imath)}\nolimits}
	\newcommand{\utMHg}{\operatorname{\ch(Q,\btau)}\nolimits}
	\newcommand{\tMHg}{\operatorname{\widetilde{\ch}(Q,\btau)}\nolimits}
	\newcommand{\tMHk}{{\widetilde{\ch}(\K Q,\btau)}}
	\newcommand{\rMHg}{\operatorname{\ch_{\rm{red}}(Q,\btau)}\nolimits}
	
	\newcommand{\rMHd}{\operatorname{\cm\ch_{\rm{red}}(\Lambda^\imath)_{\bvsd}}\nolimits}
	\newcommand{\tMHd}{\operatorname{\cs\cd\widetilde{\ch}(\Lambda^\imath)_{\bvsd}}\nolimits}
	
	\newcommand{\tMHl}{\cs\cd\widetilde{\ch}({\bs}_\ell\Lambda^\imath)}
	\newcommand{\rMHl}{\cm\ch_{\rm{red}}({\bs}_\ell\Lambda^\imath)_{\bvsd}}
	\newcommand{\tMHi}{\cs\cd\widetilde{\ch}({\bs}_i\Lambda^\imath)}
	\newcommand{\rMHi}{\cm\ch_{\rm{red}}({\bs}_i\Lambda^\imath)_{\bvsd}}
	\newcommand{\tMHgi}{\widetilde{\ch}({\bs}_i Q,\btau)}

	\newcommand{\utGpg}{\operatorname{\ch^{\rm Gp}(Q,\btau)}\nolimits}
	\newcommand{\tGpg}{\operatorname{\widetilde{\ch}^{\rm Gp}(Q,\btau)}\nolimits}
	\newcommand{\rGpg}{\operatorname{\ch_{red}^{\rm Gp}(Q,\btau)}\nolimits}

	\newcommand{\colim}{\operatorname{colim}\nolimits}
	\newcommand{\gldim}{\operatorname{gl.dim}\nolimits}
	\newcommand{\cone}{\operatorname{cone}\nolimits}
	\newcommand{\rep}{\operatorname{rep}\nolimits}
	\newcommand{\Ext}{\operatorname{Ext}\nolimits}
	\newcommand{\Tor}{\operatorname{Tor}\nolimits}
	\newcommand{\Hom}{\operatorname{Hom}\nolimits}
	\newcommand{\Top}{\operatorname{top}\nolimits}
	\newcommand{\Coker}{\operatorname{Coker}\nolimits}
	\newcommand{\thick}{\operatorname{thick}\nolimits}
	\newcommand{\rank}{\operatorname{rank}\nolimits}
	\newcommand{\Gproj}{\operatorname{Gproj}\nolimits}
	\newcommand{\Len}{\operatorname{Length}\nolimits}
	\newcommand{\RHom}{\operatorname{RHom}\nolimits}
	\renewcommand{\deg}{\operatorname{deg}\nolimits}
	\renewcommand{\Im}{\operatorname{Im}\nolimits}
	\newcommand{\Ker}{\operatorname{Ker}\nolimits}
	\newcommand{\Coh}{\operatorname{Coh}\nolimits}
	\newcommand{\Id}{\operatorname{Id}\nolimits}
	\newcommand{\Qcoh}{\operatorname{Qch}\nolimits}
	\newcommand{\CM}{\operatorname{CM}\nolimits}
	\newcommand{\sgn}{\operatorname{sgn}\nolimits}
	\newcommand{\Gdim}{\operatorname{G.dim}\nolimits}
	\newcommand{\fpr}{\operatorname{\mathcal{P}^{\leq1}}\nolimits}
	
	\newcommand{\For}{\operatorname{{\bf F}or}\nolimits}
	\newcommand{\coker}{\operatorname{Coker}\nolimits}
	\renewcommand{\dim}{\operatorname{dim}\nolimits}
	\newcommand{\rankv}{\operatorname{\underline{rank}}\nolimits}
	\newcommand{\dimv}{{\operatorname{\underline{dim}}\nolimits}}
	\newcommand{\diag}{{\operatorname{diag}\nolimits}}
	\newcommand{\qbinom}[2]{\begin{bmatrix} #1\\#2 \end{bmatrix} }
	
	\renewcommand{\Vec}{{\operatorname{Vec}\nolimits}}
	\newcommand{\pd}{\operatorname{proj.dim}\nolimits}
	\newcommand{\gr}{\operatorname{gr}\nolimits}
	\newcommand{\id}{\operatorname{Id}\nolimits}
	\newcommand{\Res}{\operatorname{Res}\nolimits}
	\def \tT{\widetilde{\mathcal T}}
	\def \tTL{\tT(\Lambda^\imath)}
	
	\newcommand{\mbf}{\mathbf}
	\newcommand{\mbb}{\mathbb}
	\newcommand{\mrm}{\mathrm}
	\newcommand{\cbinom}[2]{\left\{ \begin{matrix} #1\\#2 \end{matrix} \right\}}
	\newcommand{\dvev}[1]{{B_1|}_{\ev}^{{(#1)}}}
	\newcommand{\dv}[1]{{B_1|}_{\odd}^{{(#1)}}}
	\newcommand{\dvd}[1]{t_{\odd}^{{(#1)}}}
	\newcommand{\dvp}[1]{t_{\ev}^{{(#1)}}}
	\newcommand{\ev}{\bar{0}}
	\newcommand{\odd}{\bar{1}}
	\newcommand{\Iblack}{\I_{\bullet}}
	\newcommand{\wb}{w_\bullet}
	\newcommand{\Uidot}{\dot{\bold{U}}^{\imath}}
	
	\newcommand{\kk}{h}
	\newcommand{\la}{\lambda}
	\newcommand{\LR}[2]{\left\llbracket \begin{matrix} #1\\#2 \end{matrix} \right\rrbracket}
	\newcommand{\ff}{B}
	\newcommand{\pdim}{\operatorname{proj.dim}\nolimits}
	\newcommand{\idim}{\operatorname{inj.dim}\nolimits}
	\newcommand{\Gd}{\operatorname{G.dim}\nolimits}
	\newcommand{\Ind}{\operatorname{Ind}\nolimits}
	\newcommand{\add}{\operatorname{add}\nolimits}
	\newcommand{\pr}{\operatorname{pr}\nolimits}
	\newcommand{\oR}{\operatorname{R}\nolimits}
	\newcommand{\oL}{\operatorname{L}\nolimits}
	\newcommand{\ext}{{ \mathfrak{Ext}}}
	\newcommand{\Perf}{{\mathfrak Perf}}
	\def\scrP{\mathscr{P}}
	\newcommand{\bk}{{\mathbb K}}
	\newcommand{\cc}{{\mathcal C}}
	\newcommand{\gc}{{\mathcal GC}}
	\newcommand{\dg}{{\rm dg}}
	\newcommand{\ce}{{\mathcal E}}
	\newcommand{\cs}{{\mathcal S}}
	\newcommand{\cl}{{\mathcal L}}
	\newcommand{\cf}{{\mathcal F}}
	\newcommand{\cx}{{\mathcal X}}
	\newcommand{\cy}{{\mathcal Y}}
	\newcommand{\ct}{{\mathcal T}}
	\newcommand{\cu}{{\mathcal U}}
	\newcommand{\cv}{{\mathcal V}}
	\newcommand{\cn}{{\mathcal N}}
	\newcommand{\mcr}{{\mathcal R}}
	\newcommand{\ch}{{\mathcal H}}
	\newcommand{\ca}{{\mathcal A}}
	\newcommand{\cb}{{\mathcal B}}
	\newcommand{\ci}{{\I}_{\btau}}
	\newcommand{\cj}{{\mathcal J}}
	\newcommand{\cm}{{\mathcal M}}
	\newcommand{\cp}{{\mathcal P}}
	\newcommand{\cg}{{\mathcal G}}
	\newcommand{\cw}{{\mathcal W}}
	\newcommand{\co}{{\mathcal O}}
	\newcommand{\cq}{{Q^{\rm dbl}}}
	\newcommand{\cd}{{\mathcal D}}
	\newcommand{\ck}{\widetilde{\mathcal K}}
	\newcommand{\calr}{{\mathcal R}}
	\newcommand{\iLa}{\Lambda^{\imath}}
	\newcommand{\La}{\Lambda}
	\newcommand{\ol}{\overline}
	\newcommand{\ul}{\underline}
	\newcommand{\st}{[1]}
	\newcommand{\ow}{\widetilde}
	\renewcommand{\P}{\mathbf{P}}
	\newcommand{\pic}{\operatorname{Pic}\nolimits}
	\newcommand{\Spec}{\operatorname{Spec}\nolimits}
	
	\newtheorem{theorem}{Theorem}[section]
	\newtheorem{acknowledgement}[theorem]{Acknowledgement}
	\newtheorem{algorithm}[theorem]{Algorithm}
	\newtheorem{assumption}[theorem]{Assumption}
	\newtheorem{axiom}[theorem]{Axiom}
	\newtheorem{case}[theorem]{Case}
	\newtheorem{claim}[theorem]{Claim}
	\newtheorem{conclusion}[theorem]{Conclusion}
	\newtheorem{condition}[theorem]{Condition}
	\newtheorem{conjecture}[theorem]{Conjecture}
	\newtheorem{construction}[theorem]{Construction}
	\newtheorem{corollary}[theorem]{Corollary}
	\newtheorem{criterion}[theorem]{Criterion}
	\newtheorem{definition}[theorem]{Definition}
	\newtheorem{example}[theorem]{Example}
	\newtheorem{exercise}[theorem]{Exercise}
	\newtheorem{lemma}[theorem]{Lemma}
	\newtheorem{notation}[theorem]{Notation}
	\newtheorem{problem}[theorem]{Problem}
	\newtheorem{proposition}[theorem]{Proposition}
	\newtheorem{solution}[theorem]{Solution}
	\newtheorem{summary}[theorem]{Summary}
	\numberwithin{equation}{section}
	
	\theoremstyle{remark}
	\newtheorem{remark}[theorem]{Remark}
	\newcommand{\Pd}{\pi_*}
	\def \bvs{{\boldsymbol{\varsigma}}}
	\def \bvsd{{\boldsymbol{\varsigma}_{\diamond}}}
	\def \btau{{{\tau}}}

	\def \bp{{\mathbf p}}
	\def \bq{{\bm q}}
	\def \bv{{v}}
	\def \bs{{\bm s}}
	
	\def \bfK{{\mathbf K}}
	
	\newcommand{\tCMHg}{\cc\widetilde{\ch}(Q,\btau)}
	\newcommand{\bfv}{\mathbf{v}}
	\def \bA{{\mathbf A}}
	\def \ba{{\mathbf a}}
	\def \bL{{\mathbf L}}
	\def \bF{{\mathbf F}}
	\def \bS{{\mathbf S}}
	\def \bC{{\mathbf C}}
	\def \bU{{\mathbf U}}
	\def \bc{{\mathbf c}}
	\def \fpi{\mathfrak{P}^\imath}
	\def \Ni{N^\imath}
	\def \fp{\mathfrak{P}}
	\def \fg{\mathfrak{g}}
	\def \fk{\fg^\theta}  
	
	\def \fn{\mathfrak{n}}
	\def \fh{\mathfrak{h}}
	\def \fu{\mathfrak{u}}
	\def \fv{\mathfrak{v}}
	\def \fa{\mathfrak{a}}
	\def \fq{\mathfrak{q}}
	\def \Z{{\Bbb Z}}
	\def \F{{\Bbb F}}
	\def \D{{\Bbb D}}
	\def \C{{\Bbb C}}
	\def \N{{\Bbb N}}
	\def \Q{{\Bbb Q}}
	\def \G{{\Bbb G}}
	\def \P{{\Bbb P}}
	\def \K{{\mathbf k}}
	\def \bK{{\Bbb K}}
	
	\def \E{{\Bbb E}}
	\def \A{{\Bbb A}}
	\def \L{{\Bbb L}}
	\def \I{{\Bbb I}}
	\def \BH{{\Bbb H}}
	\def \T{{\Bbb T}}
	\newcommand{\TT}{\operatorname{\texttt{\rm T}}\nolimits}
	\newcommand {\lu}[1]{\textcolor{red}{$\clubsuit$: #1}}
	
	\newcommand{\nc}{\newcommand}
	\newcommand{\browntext}[1]{\textcolor{brown}{#1}}
	\newcommand{\greentext}[1]{\textcolor{green}{#1}}
	\newcommand{\redtext}[1]{\textcolor{red}{#1}}
	\newcommand{\bluetext}[1]{\textcolor{blue}{#1}}
	\newcommand{\brown}[1]{\browntext{ #1}}
	\newcommand{\green}[1]{\greentext{ #1}}
	\newcommand{\red}[1]{\redtext{ #1}}
	\newcommand{\blue}[1]{\bluetext{ #1}}
	
	\newcommand{\wtodo}{\todo[inline,color=orange!20, caption={}]}
	\newcommand{\lutodo}{\todo[inline,color=green!20, caption={}]}

	\title[Hall algebras and quantum symmetric pairs of Kac-Moody type II]{Hall algebras and quantum symmetric pairs of Kac-Moody type II}
	
	\author[Ming Lu]{Ming Lu}
	\address{Department of Mathematics, Sichuan University, Chengdu 610064, P.R.China}
	\email{luming@scu.edu.cn}

	\author[Runze Shang]{Runze Shang}
	\address{Department of Mathematics, Sichuan University, Chengdu 610064, P.R.China}
	\email{shangrunze@stu.scu.edu.cn }

	
	\begin{abstract}
		{We extend the $\imath$Hall algebra realization of $\imath$quantum groups arising from quantum symmetric pairs, 
			which establishes an injective homomorphism from the universal $\imath$quantum group of Kac-Moody type to the $\imath$Hall algebra associated to an arbitrary $\imath$quiver (not necessarily virtually acyclic). This generalizes Lu-Wang's result.}  
	\end{abstract}
		\maketitle
	\setcounter{tocdepth}{1}
	\tableofcontents
	
	\section{Introduction}
	$\imath$Quantum groups $\Ui$ (also called quantum symmetric pair coideal subalgebras) were introduced by G. Letzter \cite{Let99,Let02}, which are right coideal subalgebras of Drinfeld-Jimbo's quantum groups $\U$; see \cite{Ko14} for an extension to Kac-Moody type. A striking breakthrough of the $\imath$quantum groups is the discovery of canonical basis by H. Bao and W. Wang \cite{BW18a}. Inspired by the Hall algebra realization, the first author and Wang  \cite{LW19a} defined the universal $\imath$quantum groups $\tUi$, such that $\Ui$ is a quotient algebra of $\tUi$ by a certain ideal generated by central elements.
	
	A Serre presentation of {\em quasi-split} $\imath$quantum groups $\Ui$ and $\tUi$ of Kac-Moody type is more complicated than a Serre presentation (which is the definition) of a quantum group, and it was recently completed in full generality by the first author joint with X.~Chen and W.~Wang \cite{CLW18}; a complete presentation of $\Ui$ in finite type was already given earlier by Letzter \cite{Let02}. 
	A crucial relation, known as the $\imath$Serre relation, in the final presentation for $\Ui$ and $\tUi$, involves the $\imath$divided powers which arise from the theory of canonical basis for quantum symmetric pairs \cite{BW18a, BeW18}. The $\imath$divided powers come in 2 forms, depending on a parity.
	
	In \cite{LW19a,LW20a}, the first author and Wang  defined a new class of quiver algebras $\Lambda^\imath$ associated to $\imath$quivers $(Q,\btau)$, and then formulated the $\imath$Hall algebra of $\Lambda^\imath$ over a finite field $\bfk=\F_q$ in the framework of semi-derived Ringel-Hall algebras of $1$-Gorenstein algebras \cite{LP16,Lu19}. This new form of Hall algebras was motivated by the construction of Bridgeland's Hall algebra of complexes \cite{Br} and Gorsky's semi-derived Hall algebras \cite{Gor2} (which were in turn built on \cite{Rin90,Gr95}). 
	
	The $\imath$Hall algebras of $\imath$quiver algebras were conjectured to provide a realization of the universal $\imath$quantum groups arising from quasi-split quantum symmetric pairs of Kac-Moody type. This conjecture was established for finite type in \cite{LW19a}, and for Kac-Moody type by using the so-called {\em virtually acyclic} $\imath$quivers in \cite{LW20a}. 
	
	The goal of this paper is to establish the conjecture completely by using $\imath$Hall algebras of {\em arbitrary} $\imath$quivers.
	Based on \cite{LW20a}, the key is to reappear $\imath$Serre relations \eqref{relation6} in $\imath$Hall algebras for arbitrary $\imath$quivers, which is settled in Section \ref{sec:iSerre-iHall}; see Theorem \ref{thm: iserre relation}. With the help of Theorem \ref{thm: iserre relation} and using the same proof of \cite[Theorem 9.6]{LW20a}, we obtain the algebra embedding $\widetilde{\Psi}: \tUi_{|v={\sqq}}{\rightarrow} \tMHk$ in Theorem \ref{thm:Ui=iHall}.

	The paper is organized as follows. 
	We review the materials for $\imath$quiver algebras, $\imath$Hall algebras and $\imath$quantum groups in Sections \ref{sec:iquiver}--\ref{sec:braid}, to make this paper self-contained. In Section \ref{sec:iSerre-iHall}, we reappear the $\imath$Serre relations in $\imath$Hall algebras. In Section \ref{subsec:iQGiH},  we establish the algebra embedding $\widetilde{\Psi}: \tUi_{|v={\sqq}}{\rightarrow} \tMHk$ and a reduced variant.

	\noindent{\bf Acknowledgments.}
	The authors thank Liangang Peng for his continuing encouragement. ML thanks Weiqiang Wang for guiding him to study  $\imath$quantum groups and their categorical realizations, and also his continuing encouragement. ML is partially supported by the National Natural Science Foundation of China (No. 12171333).  The authors thank the referee for his/her
	careful reading and helpful suggestions.

	\section{$\imath$Quiver algebras and $\imath$Hall algebras}
	\label{sec:iquiver}
	In this section, we review the materials on $\imath$quiver algebras and $\imath$Hall algebras, following \cite{LW19a,LW20a}; also see \cite{LP16,Lu19}.

	Let $\K$ be a field.
	For a quiver algebra $A=\K Q/I$ (not necessarily finite-dimensional), we always identify left $A$-modules with representations of $Q=(Q_0,Q_1)$ satisfying relations in $I$. A representation $V=(V_i,V(\alpha))_{i\in Q_0,\alpha\in Q_1}$ of $A$ is called {\em nilpotent} if for each oriented cycle $\alpha_m\cdots\alpha_1$ at a vertex $i$, the $\K$-linear map $V(\alpha_m)\cdots V(\alpha_1):V_i\rightarrow V_i$ is nilpotent. We denote by $\mod(A)$ the category of finite-dimensional nilpotent $A$-modules.
	
	For any $A$-module $M$, we denote by $\top M$ the top of $M$, and by $\add M$ the subcategory formed by direct summands of finite direct sums of
	copies of $M$.
	
	\subsection{$\imath$Quiver algebras}

	We recall the $\imath$quiver algebras from \cite[\S2]{LW20a} and  \cite[\S2]{LW19a}.
	
	Let $\K$ be a field.
	Let $Q=(Q_0,Q_1)$ be a quiver (not necessarily acyclic). Throughout the paper, we shall identify
	\[
	\I=Q_0.
	\]
	An {\em involution} of $Q$ is defined to be an automorphism $\btau$ of the quiver $Q$ such that $\btau^2=\Id$. In particular, we allow the {\em trivial} involution $\Id:Q\rightarrow Q$. An involution $\btau$ of $Q$ induces an involution of the path algebra $\K Q$, again denoted by $\btau$.
	A quiver together with a specified involution $\btau$, $(Q, \btau)$, will be called an {\em $\imath$quiver}.
	
	Let $R_1$ denote the truncated polynomial algebra $\K[\varepsilon]/(\varepsilon^2)$.
	Let $R_2$ denote the radical square zero of the path algebra of $\xymatrix{1 \ar@<0.5ex>[r]^{\varepsilon} & 1' \ar@<0.5ex>[l]^{\varepsilon'}}$, i.e., $\varepsilon' \varepsilon =0 =\varepsilon\varepsilon '$. Define a $\K$-algebra
	\begin{equation}
		\label{eq:La}
		\Lambda=\K Q\otimes_\K R_2.
	\end{equation}
	
	Associated to the quiver $Q$, the {\em double framed quiver} $Q^\sharp$
	is the quiver such that
	\begin{itemize}
		\item the vertex set of $Q^{\sharp}$ consists of 2 copies of the vertex set $Q_0$, $\{i,i'\mid i\in Q_0\}$;
		\item the arrow set of $Q^{\sharp}$ is
		\[
		\{\alpha: i\rightarrow j,\alpha': i'\rightarrow j'\mid(\alpha:i\rightarrow j)\in Q_1\}\cup\{ \varepsilon_i: i\rightarrow i' ,\varepsilon'_i: i'\rightarrow i\mid i\in Q_0 \}.
		\]
	\end{itemize}
	Note $Q^\sharp$ admits a natural involution, $\swa$.
	The involution $\btau$ of a quiver $Q$ induces an involution ${\btau}^{\sharp}$ of $Q^{\sharp}$ which is basically the composition of $\swa$ and $\btau$ (on the two copies of subquivers $Q$ and $Q'$ of $Q^\sharp$). The algebra $\Lambda$ can be realized in terms of the quiver $Q^{\sharp}$ and a certain ideal $I^{\sharp}$
	so that $\Lambda\cong \K Q^{\sharp} \big/ I^{\sharp}$.
	
	By definition, ${\btau}^{\sharp}$ on $Q^\sharp$ preserves $I^\sharp$ and hence induces an involution ${\btau}^{\sharp}$ on the algebra $\Lambda$. The {\rm $\imath$quiver algebra} of $(Q, \btau)$ is the fixed point subalgebra of $\Lambda$ under ${\btau}^{\sharp}$,
	\begin{equation}
		\label{eq:iLa}
		\iLa
		= \{x\in \Lambda\mid {\btau}^{\sharp}(x) =x\}.
	\end{equation}
	The algebra $\iLa$ can be described in terms of a certain quiver $\ov Q$ and its ideal $\ov{I}$ so that $\iLa \cong \K \ov{Q} / \ov{I}$; see \cite[Proposition 2.6]{LW19a}.
	We recall $\ov{Q}$ and $\ov{I}$ as follows:
	\begin{itemize}
		\item $\ov{Q}$ is constructed from $Q$ by adding a loop $\varepsilon_i$ at the vertex $i\in Q_0$ if $i=\btau i$, and adding an arrow $\varepsilon_i: i\rightarrow \btau i$ for each $i\in Q_0$ if $i\neq \btau i$;
		\item $\ov{I}$ is generated by
		\begin{itemize}
			\item (Nilpotent relations) $\varepsilon_{i}\varepsilon_{\btau i}$ for any $i\in\I$;
			\item (Commutative relations) $\varepsilon_i\alpha-\btau(\alpha)\varepsilon_j$ for any arrow $\alpha:j\rightarrow i$ in $Q_1$.
		\end{itemize}
	\end{itemize}
	Moreover, it follows by \cite[Proposition 2.2]{LW20a} that $\Lambda^{\imath}$ is a $1$-Gorenstein algebra.

	By \cite[Corollary 2.12]{LW19a}, $\K Q$ is naturally a subalgebra and also a quotient algebra of $\Lambda^\imath$.
	Viewing $\K Q$ as a subalgebra of $\Lambda^{\imath}$, we have a restriction functor
	\[
	\res: \mod (\Lambda^{\imath})\longrightarrow \mod (\K Q).
	\]
	Viewing $\K Q$ as a quotient algebra of $\Lambda^{\imath}$, we obtain a pullback functor
	\begin{equation}\label{eqn:rigt adjoint}
		\iota:\mod(\K Q)\longrightarrow\mod(\Lambda^{\imath}).
	\end{equation}

	For each $i\in Q_0$, define a $k$-algebra
	\begin{align}\label{dfn:Hi}
		\BH _i:=\left\{ \begin{array}{cc}  \K[\varepsilon_i]/(\varepsilon_i^2) & \text{ if }i=\btau i,
			\\
			\K(\xymatrix{i \ar@<0.5ex>[r]^{\varepsilon_i} & \btau i \ar@<0.5ex>[l]^{\varepsilon_{\btau i}}})/( \varepsilon_i\varepsilon_{\btau i},\varepsilon_{\btau i}\varepsilon_i)  &\text{ if } \btau i \neq i .\end{array}\right.
	\end{align}
	Note that $\BH _i=\BH _{\btau i}$ for any $i\in Q_0$.
	Choose one representative for each $\btau$-orbit on $\I$, and let
	\begin{align}
		\label{eq:ci}
		\ci = \{ \text{the chosen representatives of $\btau$-orbits in $\I$} \}.
	\end{align}
	
	Define the following subalgebra of $\Lambda^{\imath}$:
	\begin{equation}  \label{eq:H}
		\BH =\bigoplus_{i\in \ci }\BH _i.
	\end{equation}
	Note that $\BH $ is a radical square zero selfinjective algebra. Denote by
	\begin{align}
		\res_\BH :\mod(\iLa)\longrightarrow \mod(\BH )
	\end{align}
	the natural restriction functor.
	On the other hand, as $\BH $ is a quotient algebra of $\iLa$ (cf. \cite[proof of Proposition~ 2.15]{LW19a}), every $\BH $-module can be viewed as a $\iLa$-module.
	
	Recall the algebra $\BH _i$ for $i \in \ci$ from \eqref{dfn:Hi}. For $i\in Q_0 =\I$, define the indecomposable module over $\BH _i$ (if $i\in \ci$) or over $\BH_{\btau i}$ (if $i\not \in \ci$)
	\begin{align}
		\label{eq:E}
		\E_i =\begin{cases}
			\bfk[\varepsilon_i]/(\varepsilon_i^2), & \text{ if }i=\btau i;
			\\
			\xymatrix{\K\ar@<0.5ex>[r]^1 & \K\ar@<0.5ex>[l]^0} \text{ on the quiver } \xymatrix{i\ar@<0.5ex>[r]^{\varepsilon_i} & \btau i\ar@<0.5ex>[l]^{\varepsilon_{\btau i}} }, & \text{ if } i\neq \btau i.
		\end{cases}
	\end{align}
	Then $\E_i$, for $i\in Q_0$, can be viewed as a $\iLa$-module and will be called a {\em generalized simple} $\iLa$-module.
	
	Let $\cp^{<\infty}(\Lambda^\imath)$ be the subcategory of $\mod(\Lambda^\imath)$ formed by modules of finite projective dimensions.
	Let $\cp^{\leq d}(\Lambda^\imath)$ be the subcategory of $\mod(\Lambda^\imath)$ which consists of $\Lambda^\imath$-modules of projective dimension less than or equal to $d$, for $d\in\N$.
	Then $\cp^{<\infty}(\Lambda^\imath)=\cp^{\leq1}(\Lambda^\imath)$, and $\E_i\in\cp^{<\infty}(\Lambda^\imath)$ for any $i\in \I$; see \cite[Lemma 2.3]{LW20a}.
	
	Following \cite{LW19a}, we can define the Euler forms $\langle K,M\rangle =\langle K,M\rangle_{\Lambda^\imath}$ and $\langle M,K\rangle =\langle M,K\rangle_{\Lambda^\imath}$ for any $K\in\cp^{\leq1}(\Lambda^\imath)$, $M\in\mod(\Lambda^\imath)$. These forms descend to bilinear Euler forms on the Grothendieck groups:
	\begin{align*}
		\langle\cdot,\cdot\rangle: K_0(\cp^{\leq 1}(\Lambda^\imath))\times K_0(\mod(\Lambda^\imath))\longrightarrow \Z,
		\\
		\langle\cdot,\cdot\rangle: K_0(\mod(\Lambda^\imath))\times K_0(\cp^{\leq 1}(\Lambda^\imath))\longrightarrow \Z.
	\end{align*}
	
	Denote by $\langle\cdot,\cdot\rangle_Q$ the Euler form of $\K Q$. Denote by $S_i$ the simple $\K Q$-module (respectively, $\Lambda^{\imath}$-module) corresponding to vertex $i\in Q_0$ (respectively, $i\in\ov{Q}_0$).
	These 2 Euler forms are related via the restriction functor $\res:\mod(\Lambda^\imath)\rightarrow \mod(\K Q)$ as follows.
	
	\begin{lemma}
		[\text{\cite[Lemma 3.1]{LW20a}}]
		\label{lemma compatible of Euler form}
		We have
		\begin{enumerate}
			\item[(1)] $\langle K,M\rangle =\langle \res_{\BH}(K),M\rangle$, $\langle M,K\rangle =\langle M,\res_{\BH}(K)\rangle$, $\forall M\in \mod(\Lambda^\imath), K\in\cp^{\leq1}(\Lambda^\imath)$;
			\item[(2)]
			$\langle \E_i, M\rangle = \langle S_i,\res (M) \rangle_Q$, $\langle M,\E_i\rangle =\langle \res(M), S_{\btau i} \rangle_Q$, $\forall i\in Q_0$, $M\in\mod(\Lambda^{\imath})$;
			\item[(3)] $\langle M,N\rangle=\frac{1}{2}\langle \res(M),\res(N)\rangle_Q$, $\forall M,N\in\cp^{\leq 1}(\Lambda^\imath)$.
		\end{enumerate}
	\end{lemma}

	\subsection{$\imath$Hall algebras}
	
	In this subsection we consider $\bfk=\F_q$, and set
	\[
	\sqq=\sqrt{q}.
	\]
	
	Generalizing \cite{LP16}, the first author defined a (twisted) semi-derived Hall algebra for a 1-Gorenstein algebra \cite{Lu19}. The $\imath$Hall algebra $\tMHk$ for $\imath$quiver $(Q,\btau)$ is by definition the twisted semi-derived Hall algebra for the module category of the $\imath$quiver algebra $\iLa$; see \cite{LW19a,LW20a}. We recall it here briefly.
	
	Let $\ch(\Lambda^\imath)$ be the Ringel-Hall algebra of $\Lambda^\imath$, i.e.,
	$$\ch(\Lambda^\imath)=\bigoplus_{[M]\in\Iso(\mod(\Lambda^\imath))}\Q(\sqq)[M],$$
	with the multiplication defined by (see \cite{Br})
	\[
	[M]\diamond [N]=\sum_{[L]\in\Iso(\mod(\Lambda^\imath))}\frac{|\Ext^1(M,N)_L|}{|\Hom(M,N)|}[L].
	\]
	
	For any three objects $X,Y,Z$, let
	\begin{align}
		F_{XY}^Z &= \big|\{L\subseteq Z, L \cong Y\text{ and }Z/L\cong X\} \big|
		\notag \\
		&= \frac{\big|\Ext^1(X,Y)_Z|}{|\Hom(X,Y)\big|} \cdot \frac{|\aut(Z)|}{|\aut(X)| |\aut(Y)|}
		\qquad \text{(Riedtman-Peng formula)}.
		\label{Ried-P}
	\end{align}
	
	Define $I$ to be the two-sided ideal of $\ch(\Lambda^\imath)$ generated by
	\begin{align}
		\label{eq:ideal}
		&\{[K]-[K'] \mid \res_\BH(K)\cong\res_\BH(K'),  K,K'\in\cp^{<\infty}(\Lambda^\imath)\} \bigcup
		\\\notag
		&\{[L]-[K\oplus M]\mid \exists \text{ exact sequence } 0 \longrightarrow K \longrightarrow L \longrightarrow M \longrightarrow 0, K\in\cp^{<\infty}(\Lambda^\imath)\}.
	\end{align}
	
	Consider the following multiplicatively closed subset $\cs$ of $\ch(\Lambda^\imath)/I$:
	\begin{equation}
		\label{eq:Sca}
		\cs = \{ a[K] \in \ch(\Lambda^\imath)/I \mid a\in \Q(\sqq)^\times, K\in \cp^{<\infty}(\Lambda^\imath)\}.
	\end{equation}
	
	The semi-derived Hall algebra of $\Lambda^\imath$ \cite{Lu19} is defined to be the localization
	$$\cs\cd\ch(\Lambda^\imath):= (\ch(\Lambda^\imath)/I)[\cs^{-1}].$$
	We define the $\imath$Hall algebra (i.e., a twisted semi-derived Hall algebra) $\tMHk$ \cite[\S4.4]{LW19a} to be the $\Q(\sqq)$-algebra on the same vector space as $\utMH$ but with twisted multiplication given by
	\begin{align}
		\label{eqn:twsited multiplication}
		[M]* [N] =\sqq^{\langle \res(M),\res(N)\rangle_Q} [M]\diamond[N].
	\end{align}

	\section{Quantum groups and $\imath$quantum groups}
	\label{sec:braid}
	
	In this section, we review the quantum groups and $\imath$quantum groups following \cite{LW19a,LW20a}; also see \cite{Lus93,Let99,Let02,Ko14,BW18a}.

	\subsection{Quantum groups}
	\label{subsec:QG}
	
	Let $Q$ be a quiver (without loops) with vertex set $Q_0= \I$.
	Let $n_{ij}$ be the number of edges connecting vertices $i$ and $j$. Let $C=(c_{ij})_{i,j \in \I}$ be the symmetric generalized Cartan matrix of the underlying graph of $Q$, defined by $c_{ij}=2\delta_{ij}-n_{ij}.$ Let $\fg$ be the corresponding Kac-Moody Lie algebra. 

	Let $v$ be an indeterminate.
	Define the quantum integers, quantum (double) factorials, and quantum binomial coefficients, for $r \in \N$ and $m \in \Z$,
	\begin{align}  \label{eq:binom}
		\begin{split}
			[m]= [m]_v =\frac{\bv^m-\bv^{-m}}{\bv-\bv^{-1}},
			&\qquad
			[r]^{!}= [r]_v^! =\prod_{i=1}^r [i]_v,
			\\
			[2r]^{!!} = [2r]_v^{!!} =\prod_{i=1}^r [2i]_v,
			& \qquad
			\qbinom{m}{r} =\frac{[m][m-1]\ldots [m-r+1]}{[r]^!}.
		\end{split}
	\end{align}
	Write $[A, B]=AB-BA$. Then $\tU = \tU_\bv(\fg)$ is defined to be the $\Q(\bv)$-algebra generated by $E_i,F_i, \tK_i,\tK_i'$, $i\in \I$, where $\tK_i, \tK_i'$ are invertible, subject to the following relations:
	\begin{align}
		[E_i,F_j]= \delta_{ij} \frac{\tK_i-\tK_i'}{\bv-\bv^{-1}},  &\qquad [\tK_i,\tK_j]=[\tK_i,\tK_j']  =[\tK_i',\tK_j']=0,
		\label{eq:KK}
		\\
		\tK_i E_j=\bv^{c_{ij}} E_j \tK_i, & \qquad \tK_i F_j=\bv^{-c_{ij}} F_j \tK_i,
		\label{eq:EK}
		\\
		\tK_i' E_j=\bv^{-c_{ij}} E_j \tK_i', & \qquad \tK_i' F_j=\bv^{c_{ij}} F_j \tK_i',
		\label{eq:K2}
	\end{align}
	and the quantum Serre relations, for $i\neq j \in \I$,
	\begin{align}
		& \sum_{r=0}^{1-c_{ij}} (-1)^r \left[ \begin{array}{c} 1-c_{ij} \\r \end{array} \right]  E_i^r E_j  E_i^{1-c_{ij}-r}=0,
		\label{eq:serre1} \\
		& \sum_{r=0}^{1-c_{ij}} (-1)^r \left[ \begin{array}{c} 1-c_{ij} \\r \end{array} \right]  F_i^r F_j  F_i^{1-c_{ij}-r}=0.
		\label{eq:serre2}
	\end{align}
	Note that $\tK_i \tK_i'$ are central in $\tU$ for all $i$.
	The comultiplication $\Delta: \widetilde{\U} \rightarrow \widetilde{\U} \otimes \widetilde{\U}$ is given by
	\begin{align}  \label{eq:Delta}
		\begin{split}
			\Delta(E_i)  = E_i \otimes 1 + \tK_i \otimes E_i, & \quad \Delta(F_i) = 1 \otimes F_i + F_i \otimes \tK_{i}', \\
			\Delta(\tK_{i}) = \tK_{i} \otimes \tK_{i}, & \quad \Delta(\tK_{i}') = \tK_{i}' \otimes \tK_{i}'.
		\end{split}
	\end{align}
	
	Analogously as for $\tU$, the quantum group $\bU$ is defined to be the $\Q(v)$-algebra generated by $E_i,F_i, K_i, K_i^{-1}$, $i\in \I$, subject to the  relations modified from \eqref{eq:KK}--\eqref{eq:serre2} with $\tK_i$ and $\tK_i'$ replaced by $K_i$ and $K_i^{-1}$, respectively. The comultiplication $\Delta$ is obtained by modifying \eqref{eq:Delta} with $\tK_i$ and $\tK_i'$ replaced by $K_i$ and $K_i^{-1}$, respectively (cf. \cite{Lus93}; beware that our $K_i$ has a different meaning from $K_i \in \U$ therein).

	\subsection{$\imath$Quantum groups}
	\label{subsec:iQG}
	
	For a  (generalized) Cartan matrix $C=(c_{ij})$, let $\Aut(C)$ be the group of all permutations $\btau$ of the set $\I$ such that $c_{ij}=c_{\btau i,\btau j}$. An element $\btau\in\Aut(C)$ is called an \emph{involution} if $\btau^2=\Id$.
	
	Let $\btau$ be an involution in $\Aut(C)$. We define $\widetilde{\bU}^\imath$ to be the $\Q(v)$-subalgebra of $\tU$ generated by
	\begin{equation}
		\label{eq:Bi}
		B_i= F_i +  E_{\btau i} \tK_i',
		\qquad \tk_i = \tK_i \tK_{\btau i}', \quad \forall i \in \I.
	\end{equation}
	Let $\tU^{\imath 0}$ be the $\Q(v)$-subalgebra of $\tUi$ generated by $\tk_i$, for $i\in \I$.
	By \cite[Lemma 6.1]{LW19a}, the elements $\tk_i$ (for $i= \btau i$) and $\tk_i \tk_{\btau i}$  (for $i\neq \btau i$) are central in $\tUi$.
	
	Let $\bvs=(\vs_i)\in  (\Q(\bv)^\times)^{\I}$ be such that $\vs_i=\vs_{\btau i}$ for all $i$. 
	Let $\Ui:=\Ui_{\bvs}$ be the $\Q(v)$-subalgebra of $\bU$ generated by
	\[
	B_i= F_i+\vs_i E_{\btau i}K_i^{-1},
	\quad
	k_j= K_jK_{\btau j}^{-1},
	\qquad  \forall i \in \I, j \in \I\backslash\ci.
	\]
	It is known \cite{Let99, Ko14} that $\bU^\imath$ is a right coideal subalgebra of $\bU$ in the sense that $\Delta: \Ui \rightarrow \Ui\otimes \U$; and $(\bU,\Ui)$ is called a \emph{quantum symmetric pair} ({\em QSP} for short), as they specialize at $v=1$ to $(U(\fg), U(\fg^{\theta}))$, where $\theta=\omega \circ \btau$, $\omega$ is the  Chevalley involution, and $\btau$ is understood here as an automorphism of $\fg$.
	
	The algebras $\Ui_{\bvs}$, for $\bvs \in  (\Q(\bv)^\times)^{\I}$, are obtained from $\tUi$ by central reductions.
	
	\begin{proposition} [\text{\cite[Proposition 6.2]{LW19a}}]
		(1) The algebra $\Ui$ is isomorphic to the quotient of $\tUi$ by the ideal generated by
		\begin{align}   \label{eq:parameters}
			\tk_i - \vs_i \; (\text{for } i =\btau i),
			\qquad  \tk_i \tk_{\btau i} - \vs_i \vs_{\btau i}  \;(\text{for } i \neq \btau i).
		\end{align}
		The isomorphism is given by sending $B_i \mapsto B_i, k_j \mapsto \vs_{\btau j}^{-1} \tk_j, k_j^{-1} \mapsto \vs_{j}^{-1} \tk_{\btau j}, \forall i\in \I, j\in \I\backslash\ci$.
		
		(2) The algebra $\widetilde{\bU}^\imath$ is a right coideal subalgebra of $\widetilde{\bU}$; that is, $(\widetilde{\bU}, \widetilde{\bU}^\imath)$ forms a QSP.
	\end{proposition}
	
	We shall refer to $\tUi$ and $\Ui$ as {\em (quasi-split) $\imath${}quantum groups}; they are called {\em split} if $\btau =\Id$.

	For $i\in \I$ with $\btau i= i$, generalizing the constructions in \cite{BW18a, BeW18}, we define the {\em $\imath${}divided powers} of $B_i$ to be (also see  \cite{CLW21})
	\begin{eqnarray}
		&&\ff_{i,\odd}^{(m)}=\frac{1}{[m]^!}\left\{ \begin{array}{ccccc} B_i\prod_{s=1}^k (B_i^2-\tk_i[2s-1]^2 ) & \text{if }m=2k+1,\\
			\prod_{s=1}^k (B_i^2-\tk_i[2s-1]^2) &\text{if }m=2k; \end{array}\right.
		\label{eq:iDPodd} \\
		&&\ff_{i,\ev}^{(m)}= \frac{1}{[m]^!}\left\{ \begin{array}{ccccc} B_i\prod_{s=1}^k (B_i^2-\tk_i[2s]^2 ) & \text{if }m=2k+1,\\
			\prod_{s=1}^{k} (B_i^2-\tk_i[2s-2]^2) &\text{if }m=2k. \end{array}\right.
		\label{eq:iDPev}
	\end{eqnarray}
	On the other hand, for $i\in \I$ with $i \neq \tau i$, we define the divided powers as in the quantum group setting: for $m\in \N$,
	\begin{align}
		\label{eq:DP}
		\ff_{i}^{(m)}= \frac{\ff_i^m}{[m]^!}.
	\end{align}
	
	Denote
	\[
	(a;x)_0=1, \qquad (a;x)_n =(1-a)(1-ax)  \cdots (1-ax^{n-1}), \quad   n\ge 1.
	\]
	
	We have the following {\em Serre presentation} of $\tUi$, with $\ov{p}_i\in \Z_2$ fixed for each $i\in \I$.
	
	\begin{proposition} [\text{\cite[Theorem 4.2]{LW20a}; also cf. \cite[Theorem 2]{CLW18}}]
		\label{prop:Serre}
		The $\Q(v)$-algebra $\tUi$ has a presentation with generators $B_i$, $\tk_i$ $(i\in \I)$ and the relations \eqref{relation1}--\eqref{relation6} below: for $\l \in \I$, and $i\neq j \in \I$,
		\begin{align}
			\tk_i \tk_l =\tk_l \tk_i,
			\quad
			\tk_i B_l & = v^{c_{\btau i,l} -c_{i l}} B_l \tk_i,
			\label{relation1}
			\\
			B_iB_{j}-B_jB_i =0, \quad &\text{ if }c_{ij} =0 \text{ and }\btau i\neq j,\label{relation2}
			\\
			\sum_{n=0}^{1-c_{ij}} (-1)^nB_i^{(n)}B_jB_i^{(1-c_{ij}-n)} &=0, \quad \text{ if }  j \neq \tau i\neq i, \label{relation3}
			\\
			\sum_{n=0}^{1-c_{i,\btau i}} (-1)^{n+c_{i,\btau i}}B_i^{(n)}B_{\btau i}&B_i^{(1-c_{i,\btau i}-n)} =\frac{1}{v-v^{-1}}\times
			\label{relation5}     \\
			\Big(v^{c_{i,\btau i}} (v^{-2};v^{-2})_{-c_{i,\btau i}}    &
			B_i^{(-c_{i,\btau i})} \tk_i
			-(v^{2};v^{2})_{-c_{i,\btau i}}B_i^{(-c_{i,\tau i})} \tk_{\btau i}  \Big),\quad
			\text{ if } \btau i \neq i,
			\notag \\
			\sum_{n=0}^{1-c_{ij}} (-1)^n  B_{i, \overline{p_i}}^{(n)}B_j B_{i,\overline{c_{ij}}+\overline{p}_i}&^{(1-c_{ij}-n)} =0,\quad   \text{ if }i=\btau i.
			\label{relation6}
		\end{align}
	\end{proposition}

	\section{$\imath$Serre relation in $\imath$Hall algebras}
	\label{sec:iSerre-iHall}
	
	In this section, we shall reappear the $\imath$Serre relation \eqref{relation6} in $\imath$Hall algebras.
	
	\subsection{$\imath$Divided powers in $\imath$Hall algebras}

	
	For a $\iLa$-module $M$, we shall write
	\[
	[l M] =[\underbrace{M\oplus \cdots \oplus M}_{l}],
	\qquad
	[M]^{l} = \underbrace{[M]* \cdots * [M]}_{l}.
	\]
	
	Following \cite{LW20a}, if $\tau i=i$, we define the $\imath$divided power of $[S_i]$ in $\tMH$ as follows: 
	\begin{align}
		\label{eq:idividedHallodd}
		&[S_i]_{\odd}^{(m)}:=\frac{1}{[m]_{\sqq}^!}\left\{ \begin{array}{ll} [S_i]*\prod_{j=1}^k ([S_i]^2+\sqq^{-1}(\sqq^2-1)^2[2j-1]_{\sqq}^2 [\E_i] ) & \text{if }m=2k+1,
			\\
			\prod_{j=1}^k ([S_i]^2+\sqq^{-1}(\sqq^2-1)^2[2j-1]_{\sqq}^2[\E_i]) &\text{if }m=2k; \end{array}\right.
		\\
		\label{eq:idividedHallev}
		&[S_i]_{\ev}^{(m)}:= \frac{1}{[m]_{\sqq}^!}\left\{ \begin{array}{ll} [S_i]*\prod_{j=1}^k ([S_i]^2+\sqq^{-1}(\sqq^2-1)^2[2j]_{\sqq}^2[\E_i] ) &\text{if }m=2k+1,\\
			\prod_{j=1}^{k} ([S_i]^2+\sqq^{-1}(\sqq^2-1)^2[2j-2]_{\sqq}^2[\E_i]) &\text{if }m=2k. \end{array}\right.
	\end{align}
	If $\tau i\neq i$, define 
	\begin{align}
		\label{eq:dividedpowerHall}
		[S_i]^{(m)}:= \frac{1}{[m]^!_\sqq} [S_i]^m=\sqq^{-\frac{m(m-1)}{2}}[m S_i] .
	\end{align}
	
	We recall the expansion formula of the $\imath$divided powers in terms of an $\imath$Hall basis. See \eqref{eq:binom} for notation $[2k]_\sqq^{!!}$.
	
	\begin{lemma}
		[\text{\cite[Propositions 6.4--6.5]{LW20a}}]
		\label{lem:iDPev}
		For any $m \in \N$, $\ov{p}\in\Z_2$ we have
		\begin{align}
			&[S_i]^{(m)}_{\ov{p}}=\begin{cases}
				\sum\limits_{k=0}^{\lfloor\frac{m}{2}\rfloor} \frac{\sqq^{k(k-1)-\binom{m-2k}{2}}}{[m-2k]_{\sqq}^{!}[2k]_\sqq^{!!}}  (\sqq-\sqq^{-1})^k
				[(m-2k)S_i]*[\E_i]^k, & \text{ if }\ov{m}=\ov{p};
				\label{eq:2mev}
				\\
				\\
				\sum\limits_{k=0}^{\lfloor\frac{m}{2}\rfloor}  \frac{\sqq^{k(k+1) -\binom{m-2k}{2}}}{[m-2k]_{\sqq}^{!}[2k]_\sqq^{!!}} (\sqq-\sqq^{-1})^k
				[(m-2k)S_i]*[\E_i]^k, & \text{ if }\ov{m}\neq\ov{p}.
			\end{cases}
		\end{align}
	\end{lemma}

	\subsection{$\imath$Serre relations in $\imath$Hall algebras}
	
	Consider the $\imath$quiver
	\begin{align}
		\label{diag:split2}
		Q=(\xymatrix{1\ar@<1ex>[r]|-{a}_{\cdots}  & 2\ar@<1ex>[l]|-{b}}),
		\quad
		\btau=\Id,
		\qquad \text{where } a+b=-c_{12}.
	\end{align}
	Here and blow, $\xymatrix{1\ar[r]|-a&2}$ means there are $a$ arrows from $1$ to $2$.
	The arrows $1\rightarrow 2$ are denoted by $\alpha_i$, $1\leq i\leq a$, and $2\rightarrow 1$ are denoted by $\beta_j$, $1\leq j\leq b$.
	Then the corresponding $\imath$quiver algebra $\Lambda^\imath$ has its quiver $\ov{Q}$ as
	\begin{center}\setlength{\unitlength}{0.8mm}
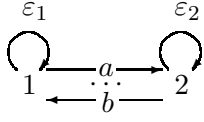
		\begin{equation*}
			\begin{picture}(50,14)(0,-8)
				\put(0,-3){$1$}
				\put(4,1){\line(1,0){8}}
				\put(12.5,0){$a$}
				\put(15.5,1){\vector(1,0){7.5}}
				\put(23,-4){\line(-1,0){7}}
				\put(12,-4){\vector(-1,0){8}}
				\put(13,-6){$b$}
				\put(11,-2.5){$\cdots$}
				\put(25,-3){$2$}
				\qbezier(-1,1)(-3,3)(-2,5.5)
				\qbezier(-2,5.5)(1,9)(4,5.5)
				\qbezier(4,5.5)(5,3)(3,1)
				\put(3.1,1.4){\vector(-1,-1){0.3}}
				\qbezier(24,1)(22,3)(23,5.5)
				\qbezier(23,5.5)(26,9)(29,5.5)
				\qbezier(29,5.5)(30,3)(28,1)
				\put(28.1,1.4){\vector(-1,-1){0.3}}
				\put(0,10){$\varepsilon_1$}
				\put(25,10){$\varepsilon_2$}
			\end{picture}
		\end{equation*}
	\end{center}

	\begin{theorem}
		\label{thm: iserre relation}
		Let $\Lambda^\imath$ be the $\imath$quiver algebra associated to the $\imath$quiver \eqref{diag:split2}.
		Then the following identity holds in $\tMHk$, for any $\overline{p}\in\Z/2$:
		\begin{align}
			\sum_{n=0}^{1+a+b} (-1)^n  [S_1]_{\overline{p}}^{(n)}*[S_2] *[S_1]_{\overline{a}+\ov{b}+\overline{p}}^{(1+a+b-n)} =0.
			\label{eqn:iserre1}
		\end{align}
	\end{theorem}
	
	The proof of Theorem \ref{thm: iserre relation} shall occupy the remainder of this section.

	\subsection{A building block}
	We recall a useful multiplication formula of $\imath$Hall algebras in the following lemma.
	\begin{lemma}[\text{\cite[Proposition 3.10]{LW20a}}]
		\label{prop:iHallmult}
		Let $(Q,\tau)$ be an $\imath$quiver with $\tau=\Id$.
		For any $A,B\in\mod(\bfk Q)\subset \mod(\Lambda^\imath)$, we have
		\begin{align}
			\label{Hallmult1}
			[A]*[B]=&
			\sum_{[L],[M],[N]\in\Iso({\mod(\bfk Q)})} \sqq^{-\langle A,B\rangle_Q}  q^{\langle N,L\rangle_Q }\frac{|\Ext^1(N, L)_{M}|}{|\Hom(N,L)| }
			\\
			\notag
			&\cdot |\{s\in\Hom(A,B)\mid \ker s\cong N, \coker s\cong L
			\}|\cdot [M]*[\bK_{\widehat{A}-\widehat{N}}]
		\end{align}
		in $\tMHk$.
	\end{lemma}

	For any $\Lambda^\imath$-module $M=(M_i, M(\alpha_j), M(\beta_k), M(\varepsilon_i))_{i=1,2;1\leq j\leq a,1\leq k\leq b}$ such that $\dim_\bfk M_2=1$, we define
	\begin{align}
		\label{eq:UW}
		U_M:= \bigcap_{1\leq i\leq a} \Ker M(\alpha_i),\qquad W_M:= \sum_{j=1}^b \Im M(\beta_j),
	\end{align}
	and let
	\begin{align}  \label{eq:uw}
		u_M:=\dim U_M,\qquad w_M:=\dim W_M.
	\end{align}
	For any $r\geq0$, define
	\begin{align}
		\cI_{r}:=&\{[M]\in\Iso(\mod(\bfk Q))\mid \widehat{M}=r\widehat{S_1}+\widehat{S_2}, W_M\subseteq U_M\},
		\\
		\cJ_r:=&\{[M]\in\Iso(\mod(\bfk Q))\mid \widehat{M}=r\widehat{S_1}+\widehat{S_2}, S_2\in\add \top(M)\}.
	\end{align}
	The following proposition is a generalization of \cite[Proposition 7.3]{LW20a}.
	\begin{proposition}
		\label{prop:SSSM}
		The following identity holds in $\tMHk$, for $s, t \ge 0$:
		\begin{align}
			[s S_1]*[S_2]*[t S_1]
			\label{eq:SSS}=&\sqq^{-tb}\sum_{r=0}^{\min\{s,t\}}\sum\limits_{[M]\in\cI_{s+t-2r}}\sqq^{\tilde{p}(a,b,r,s,t)}(\sqq-\sqq^{-1})^{s-r+t+1}
			\frac{[s]_\sqq^![t]_\sqq^!}{[r]_\sqq^!} \notag
			\\\quad&\begin{bmatrix} u_M-w_M \\ (t-r)-w_M \end{bmatrix}_\sqq\frac{1}{|\Aut(M)|}[M]*[\E_1]^r,
		\end{align}
		where
		\begin{align*}
			\tilde{p}(a,b,r,s,t)=&2(s-r)(t-r-a)-s(t-a)+2br-2r(t-r)-r^2+rs+t^2+\tbinom{t}{2}+
			\\&(s-r)^2+\tbinom{s-r}{2}+\big(u_M-(t-r)\big)\big((t-r)-w_M\big)+1.
		\end{align*}
	\end{proposition}
	
	\begin{proof}
		Using Lemma \ref{prop:iHallmult}, we compute
		\begin{align*}
			&\text{LHS}\eqref{eq:SSS}
			\\
			=&\sqq^{-tb} \sum\limits_{[U]\in\cJ_t} |\Ext^1(S_2,S_1^{\oplus{t}})_U| [S_1^{\oplus{s}}]*[U]
			\\
			=&\sqq^{-tb} \sum\limits_{[U]\in\cJ_t} |\Ext^1(S_2,S_1^{\oplus{t}})_U| \sum_{r=0}^{\min\{s,t\}} \sum\limits_{[M],[L]\in\cJ_{t-r}}\sqq^{sa-st}q^{(s-r)(t-r-a)}
			\\
			&\quad\frac{|\Ext^1(S_1^{\oplus{(s-r)}},L)_M|}{|\Hom(S_1^{\oplus{(s-r)}},L)|} \cdot \lvert\{f:S_1^{\oplus{t}}\rightarrow U\mid \ker f\cong S_1^{\oplus{(s-r)}},\coker f\cong L\}\rvert \cdot [M]*[\bbK_1]^r.
		\end{align*}
		
		In fact, for a fixed $0\leq r\leq \min\{s,t\}$, we have $\widehat{M}= (s+t-2r)\widehat{S_1}+\widehat{S_2}$.
		For any $f: S_1^{\oplus{s}}\rightarrow U$ such that $\ker f\cong S_1^{\oplus{(s-r)}}$ and $\coker f\cong L$, there exist $f_1:S_1^{\oplus s}\twoheadrightarrow S_1^{\oplus r}$ and $f_2:S_1^{\oplus r}\hookrightarrow U$ such that $\Ker f_1\cong S_1^{\oplus{(s-r)}}$ and $\coker f_2\cong L$. Moreover, $f_1$ and $f_2$ are unique up to a free action of $\Aut(S_1^{\oplus r})$. Then
		\begin{align*}
			&\lvert\{f:S_1^{\oplus{s}}\rightarrow U\mid \ker f\cong S_1^{\oplus{(s-r)}},\coker f\cong L\}\rvert
			\\
			=& \frac{|\{f_1:S_1^{\oplus s}\twoheadrightarrow S_1^{\oplus r}\mid \Ker f_1\cong
				S_1^{\oplus{(s-r)}} \}|\cdot |\{ f_2:S_1^{\oplus r}\hookrightarrow U\mid \coker f_2\cong L\}|}{|\Aut(S_1^{\oplus r})|}\\
			=&|\{f_1:S_1^{\oplus s}\twoheadrightarrow S_1^{\oplus r} \}|\cdot F_{L,S_1^{\oplus r}}^U.
		\end{align*}
		So
		\begin{align*}
			\text{LHS}\eqref{eq:SSS}
			=&\sqq^{-tb}\sum\limits_{r\geq0,[M]} \sum_{[U]\in\cJ_t,[L]\in\cJ_{t-r}}\sqq^{sa-st}q^{(s-r)(t-r-a)}\cdot|\Ext^1(S_2,S_1^{\oplus{t}})_U| \cdot
			\\& \quad |\Ext^1(S_1^{\oplus{(s-r)}},L)_M| \cdot q^{-(s-r)(t-r)}\cdot\lvert \{f_1:S_1^{\oplus{s}}\twoheadrightarrow S_1^{\oplus{r}}\}\rvert \cdot F_{L,S_1^{\oplus{r}}}^{U}\cdot [M]*[\bbK_1]^r
			\\=&\sqq^{-tb}\sum\limits_{r\geq0,[M]}\sum_{[U]\in\cJ_t,[L]\in\cJ_{t-r}} \sqq^{-sa-st+2ar}\prod_{j=0}^{r-1}(q^s-q^j)\cdot\lvert \Ext^1(S_2,S_1^{\oplus{t}})_U \rvert \cdot\\
			&\qquad  F_{L,S_1^{\oplus{r}}}^{U}\cdot [M]*[\bbK_1]^r.
		\end{align*}
		It is obvious that $F_{S_2,S_1^{\oplus t}}^U\neq0$ if and only if $[U]\in\cJ_t$, and in this case, $F_{S_2,S_1^{\oplus t}}^U=1$. By the Riedtmann-Peng formula we have
		\begin{align*}
			|\Ext^1(S_2,S_1^{\oplus{t}})_U |= \frac{|\Aut(S_2)|\cdot|\Aut(S_1^{\oplus t})|}{|\Aut(U)|} F_{S_2,S_1^{\oplus t}}^U= \frac{|\Aut(S_2)|\cdot|\Aut(S_1^{\oplus t})|}{|\Aut(U)|}.
		\end{align*}
		Furthermore, we have
		\begin{align*}
			\lvert \Ext^1(S_1^{\oplus{(s-r)}},L)_M \rvert =F_{S_1^{\oplus{(s-r)}},L}^{M}\cdot \lvert \Hom(S_1^{\oplus{(s-r)}},L)\rvert
			\frac{|\Aut{(S_1^{\oplus{(s-r)}})}|\cdot|\Aut(L)|}{|\Aut(M)|}
		\end{align*}
		by using the Riedtmann-Peng formula again.
		So
		\begin{align*}
			&\text{LHS}\eqref{eq:SSS}
			\\=&\sqq^{-tb}\sum_{r\geq0,[M]}\sum\limits_{[U]\in\cJ_t,[L]\in\cJ_{t-r}}\sqq^{-sa-st+2ar}\prod_{j=0}^{r-1}(q^s-q^j)\cdot\lvert \Ext^1(S_2,S_1^{\oplus{t}})_U \rvert\cdot
			F_{S_1^{\oplus{(s-r)}},L}^{M}\cdot 
			\\& \quad \lvert \Hom(S_1^{\oplus{(s-r)}},L)\rvert\cdot\frac{|\Aut{(S_1^{\oplus{(s-r)}})}|\cdot|\Aut(L)|}{|\Aut(M)|}\cdot F_{L,S_1^{\oplus{r}}}^{U}\cdot [M]*[\bbK_1]^r
			\\=&\sqq^{-tb}\sum_{r\geq0,[M]}\sum\limits_{[U]\in\cJ_t,[L]\in\cJ_{t-r}}\sqq^{2(s-r)(t-r-a)-s(t-a)}\prod_{j=0}^{r-1}(q^s-q^j)\prod_{j=0}^{s-r-1}(q^{s-r}-q^j)\cdot
			\\& \quad \lvert \Ext^1(S_2,S_1^{\oplus{t}})_U \rvert \cdot F_{S_1^{\oplus{(s-r)}},L}^{M}\cdot F_{L,S_1^{\oplus{r}}}^{U}\cdot \frac{|\Aut(L)|}{|\Aut(M)|} \cdot [M]*[\bbK_1]^r
			\\
			=&\sqq^{-tb}\sum_{r\geq0,[M]}\sum\limits_{[U]\in\cJ_t,[L]\in\cJ_{t-r}}\sqq^{2(s-r)(t-r-a)-s(t-a)}\prod_{j=0}^{r-1}(q^s-q^j)\prod_{j=0}^{s-r-1}(q^{s-r}-q^j)
			\\& \quad \frac{|\Aut(S_2)|\cdot|\Aut(S_1^{\oplus t})|}{|\Aut(U)|}
			F_{S_1^{\oplus{(s-r)}},L}^{M}\cdot F_{L,S_1^{\oplus{r}}}^{U}\cdot \frac{|\Aut(L)|}{|\Aut(M)|}\cdot [M]*[\bbK_1]^r.
		\end{align*}
		
		Since \begin{align*}
			F_{L,S_1^{\oplus{r}}}^{U}=\frac{|\Ext^1(L,S_1^{\oplus r})_U|}{|\Hom(L,S_1^{\oplus r})|} \frac{|\Aut(U)|}{|\Aut(S_1^{\oplus r})|\cdot|\Aut(L)|}.
		\end{align*}
		We have
		\begin{align*}
			&\text{LHS}\eqref{eq:SSS}
			\\
			=&\sqq^{-tb}\sum_{r\geq0,[M]}\sqq^{2(s-r)(t-r-a)-s(t-a)}\prod_{j=0}^{r-1}(q^s-q^j)\prod_{j=0}^{s-r-1}(q^{s-r}-q^j)
			\frac{|\Aut(S_2)|\cdot|\Aut(S_1^{\oplus t})|}{|\Aut(M)|\cdot|\Aut(S_1^{\oplus r})|}
			\\& \quad\sum\limits_{[L]\in\cJ_{t-r}} F_{S_1^{\oplus{(s-r)}},L}^{M}\cdot \sum\limits_{[U]\in\cJ_t}\frac{|\Ext^1(L,S_1^{\oplus r})_U|}{|\Hom(L,S_1^{\oplus r})|}  \cdot [M]*[\bbK_1]^r.
		\end{align*}
		
		For any $[L]\in \cJ_{t-r}$, we have $[U]\in \cJ_t$ for any $U$ such that $|\Ext^1(L,S_1^{\oplus r})_U|\neq 0$. Then
		\begin{align*}
			\sum\limits_{[U]\in\cJ_{t}}\frac{|\Ext^1(L,S_1^{\oplus r})_U|}{|\Hom(L,S_1^{\oplus r})|}=\frac{|\Ext^1(L,S_1^{\oplus r})|}{|\Hom(L,S_1^{\oplus r})|}=q^{-\langle L,S_1^{\oplus r}\rangle_Q}=q^{-\langle M,S_1^{\oplus r}\rangle_Q +\langle S_1^{\oplus(s-r)},S_1^{\oplus r}\rangle_Q}.
		\end{align*}
		So
		\begin{align*}
			&\text{LHS}\eqref{eq:SSS}
			\\
			=&\sqq^{-tb}\sum\limits_{[M],r}\sqq^{2(s-r)(t-r-a)-s(t-a)}\prod_{j=0}^{r-1}(q^s-q^j)\prod_{j=0}^{s-r-1}(q^{s-r}-q^j)
			\frac{|\Aut(S_2)|\cdot|\Aut(S_1^{\oplus t})|}{|\Aut(M)|}
			\\& \quad q^{-\langle M,S_1^{\oplus r}\rangle_Q +\langle S_1^{\oplus(s-r)},S_1^{\oplus r}\rangle_Q} \frac{1}{|\Aut(S_1^{\oplus r})|}\sum\limits_{[L]\in\cJ_{t-r}}F_{S_1^{\oplus{(s-r)}},L}^{M} \cdot [M]*[\bbK_1]^r.
		\end{align*}

		Similar to \cite[Page 38]{Sch12}, one can obtain
		
		\begin{align}
			\label{eq:coeffM}
			\sum\limits_{[L]\in\cJ_{t-r}}F_{S_1^{\oplus{(s-r)}},L}^{M}=&|\Gr(t-r-w_M,u_M-w_M)|
			\\\notag
			=&\sqq^{(u_M-(t-r))((t-r)-w_M)}\begin{bmatrix} u_M-w_M \\ (t-r)-w_M \end{bmatrix}_v.
		\end{align}
		Here $\Gr(t-r-w_M,u_M-w_M)$ is the Grassmanian of $(t-r-w_M)$-planes in $(u_M-w_M)$-space. Furthermore, \eqref{eq:coeffM} is not zero only if $W_M\subseteq U_M$. So we can assume $[M]\in\cI_{s+t-2r}$.
		
		A direct computation shows 
		\begin{align*}
			|\Aut(S_2)|=q-1,\qquad |\Aut(S_1^{\oplus t})|=\prod_{i=0}^{t-1}(q^t-q^i),\\
			\langle S_1^{\oplus(s-r)},S_1^{\oplus r}\rangle_Q=(s-r)r,\qquad \langle M,S_1^{\oplus r}\rangle_Q=r(s+t-2r)-br.
		\end{align*}
		
		Recall $q=\sqq^2$. Note that
		\begin{align*}
			\prod_{j=0}^{r-1}(q^r-q^j)&=\sqq^{r^2+\tbinom{r}{2}}(\sqq-\sqq^{-1})^r[r]_v^!,
			\\
			\prod_{j=0}^{t-1}(q^t-q^j)&=\sqq^{t^2+\tbinom{t}{2}}(\sqq-\sqq^{-1})^t[t]_v^!,
			\\
			\prod_{j=0}^{r-1}(q^s-q^j)&=\sqq^{rs+\tbinom{r}{2}}(\sqq-\sqq^{-1})^r[s]_v[s-1]_v\cdots [s-r+1]_v,
			\\
			\prod_{j=0}^{s-r-1}(q^{s-r}-q^j)&=\sqq^{(s-r)^2+\tbinom{s-r}{2}}(\sqq-\sqq^{-1})^{s-r}[s-r]_v^!.
		\end{align*}
		
		According to the above computations, the desired formula follows.
	\end{proof}
	
	
	\subsection{Proof of \eqref{eqn:iserre1}}
	

	Note that \eqref{eqn:iserre1} can be written as
	\begin{align}
		\sum_{n=0}^{1+a+b} (-1)^n  [S_1]_{\overline{0}}^{(n)}*[S_2] *[S_1]_{\overline{a}+\ov{b}}^{(1+a+b-n)} =0,
		\label{eqn:iserre10}
		\\
		\sum_{n=0}^{1+a+b} (-1)^n  [S_1]_{\overline{1}}^{(n)}*[S_2] *[S_1]_{\overline{a}+\ov{b}+\ov{1}}^{(1+a+b-n)} =0.
		\label{eqn:iserre11}
	\end{align}
	
	In the following, we give a detailed proof of \eqref{eqn:iserre10}, and then a sketch proof of \eqref{eqn:iserre11}. The proof given here is modified from \cite[\S 7.3.1]{LW20a}.
	We divide the computation of the \text{LHS} of \eqref{eqn:iserre10} into 2 cases.
	
	\underline{\textit{Case (I): n even}}. According to Proposition \ref{prop:SSSM} and Lemma \ref{lem:iDPev}, we have
	\begin{align*}
		&[S_1]_{\overline{0}}^{(n)}*[S_2] *[S_1]_{\overline{a}+\ov{b}}^{(1+a+b-n)}
		\\=&\sum_{k=0}^{\frac{n}{2}} \sum_{m=0}^{\lfloor\frac{a+b+1-n}{2}\rfloor} \frac{\sqq^{k(k-1)+m(m+1)-\tbinom{n-2k}{2}-\tbinom{1+a+b-n-2m}{2}}\cdot(\sqq-\sqq^{-1})^{k+m}}{[n-2k]_\sqq^![1+a+b-n-2m]_\sqq^![2k]_\sqq^{!!}[2m]_\sqq^{!!}}
		\\&\times[(n-2k)S_1]*[S_2]*[(1+a+b-n-2m)S_1]*[\E_1]^{k+m}
		\\=&\sum_{k=0}^{\frac{n}{2}} \sum_{m=0}^{\lfloor\frac{a+b+1-n}{2}\rfloor} \sqq^{-(1+a+b-n-2m)b} \sum_{r=0}^{\min\{n-2k,1+a+b-n-2m\}}\sum\limits_{[M]\in\cI_{1+a+b-2k-2m-2r}}
		\\&\frac{\sqq^{k(k-1)+m(m+1)-\tbinom{n-2k}{2}-\tbinom{1+a+b-n-2m}{2}}}{[n-2k]_\sqq^![1+a+b-n-2m]_\sqq^![2k]_\sqq^{!!}[2m]_\sqq^{!!}} (\sqq-\sqq^{-1})^{k+m} \sqq^{\tilde p(a,b,r,n-2k,1+a+b-n-2m)}
		\\&(\sqq-\sqq^{-1})^{2+a+b-2k-2m-r}\frac{[n-2k]_\sqq^![1+a+b-n-2m]_\sqq^!}{[r]_\sqq^!} 
		\\&\cdot\begin{bmatrix} u_M-w_M \\ 1+a+b-n-2m-r-w_M \end{bmatrix}_\sqq\cdot\frac{[M]}{|\Aut(M)|}*[\E_1]^{k+m+r}.
	\end{align*}
	
	This can be simplified to be, for $n$ even,
	\begin{align}
		&[S_1]_{\overline{0}}^{(n)}*[S_2] *[S_1]_{\overline{a}+\ov{b}}^{(1+a+b-n)}\label{eqn:iserre101}
		\\=&\sum_{k=0}^{\frac{n}{2}} \sum_{m=0}^{\lfloor\frac{a+b+1-n}{2}\rfloor} \sqq^{-(1+a+b-n-2m)b} \sum_{r=0}^{\min\{n-2k,1+a+b-n-2m\}}\sum\limits_{[M]\in\cI_{1+a+b-2k-2m-2r}} \notag
		\\&\frac{\sqq^{\tilde z}(\sqq-\sqq^{-1})^{2+a+b-k-m-r}}{[2k]_\sqq^{!!}[2m]_\sqq^{!!}[r]_\sqq^{!}}\begin{bmatrix} u_M-w_M \\ 1+a+b-n-2m-r-w_M \end{bmatrix}_\sqq
		\frac{[M]*[\E_1]^{k+m+r}}{|\Aut(M)|},\notag
	\end{align}
	where we denote
	\begin{align*}
		\tilde z=k(k-1)&+m(m+1)-\binom{n-2k}{2}-\binom{1+a+b-n-2m}{2}
		\\&+\tilde p(a,b,r,n-2k,1+a+b-n-2m).
	\end{align*}
	
	\underline{\textit{Case (II): n odd}}.
	\begin{align*}
		&[S_1]_{\overline{0}}^{(n)}*[S_2] *[S_1]_{\overline{a}+\ov{b}}^{(1+a+b-n)}
		\\=&\sum_{k=0}^{\frac{n}{2}} \sum_{m=0}^{\lfloor\frac{a+b+1-n}{2}\rfloor} \frac{\sqq^{k(k+1)+m(m-1)-\tbinom{n-2k}{2}-\tbinom{1+a+b-n-2m}{2}}\cdot(\sqq-\sqq^{-1})^{k+m}}{[n-2k]_\sqq^![1+a+b-n-2m]_\sqq^![2k]_\sqq^{!!}[2m]_\sqq^{!!}}
		\\&\times[(n-2k)S_1]*[S_2]*[(1+a+b-n-2m)S_1]*[\E_1]^{k+m}
		\\=&\sum_{k=0}^{\frac{n}{2}} \sum_{m=0}^{\lfloor\frac{a+b+1-n}{2}\rfloor} \sqq^{-(1+a+b-n-2m)b} \sum_{r=0}^{\min\{n-2k,1+a+b-n-2m\}}\sum\limits_{[M]\in\cI_{1+a+b-2k-2m-2r}}
		\\&\frac{\sqq^{k(k+1)+m(m-1)-\tbinom{n-2k}{2}-\tbinom{1+a+b-n-2m}{2}}}{[n-2k]_\sqq^![1+a+b-n-2m]_\sqq^![2k]_\sqq^{!!}[2m]_\sqq^{!!}} (\sqq-\sqq^{-1})^{k+m} \sqq^{\tilde p(a,b,r,n-2k,1+a+b-n-2m)}
		\\&(\sqq-\sqq^{-1})^{2+a+b-2k-2m-r}\frac{[n-2k]_\sqq^![1+a+b-n-2m]_\sqq^!}{[r]_\sqq^!} 
		\\&\cdot\begin{bmatrix} u_M-w_M \\ 1+a+b-n-2m-r-w_M \end{bmatrix}_\sqq\cdot\frac{[M]}{|\Aut(M)|}*[\E_1]^{k+m+r}.
	\end{align*}
	This can be simplified to be, for $n$ odd,
	\begin{align}
		&[S_1]_{\overline{0}}^{(n)}*[S_2] *[S_1]_{\overline{a}+\ov{b}}^{(1+a+b-n)}\label{eqn:iserre102}
		\\=&\sum_{k=0}^{\frac{n}{2}} \sum_{m=0}^{\lfloor\frac{a+b+1-n}{2}\rfloor} \sqq^{-(1+a+b-n-2m)b} \sum_{r=0}^{\min\{n-2k,1+a+b-n-2m\}}\sum\limits_{[M]\in\cI_{1+a+b-2k-2m-2r}} \notag
		\\&\frac{\sqq^{\tilde z+2k-2m}(\sqq-\sqq^{-1})^{2+a+b-k-m-r}}{[2k]_\sqq^{!!}[2m]_\sqq^{!!}[r]_\sqq^{!}}\begin{bmatrix} u_M-w_M \\ 1+a+b-n-2m-r-w_M \end{bmatrix}_\sqq
		\frac{[M]*[\E_1]^{k+m+r}}{|\Aut(M)|}. \notag
	\end{align}
	
	Summing up \eqref{eqn:iserre101} and \eqref{eqn:iserre102} above, we obtain
	\begin{align}
		&\sum_{n=0}^{a+b+1}(-1)^n[S_1]_{\overline{0}}^{(n)}*[S_2] *[S_1]_{\overline{a}+\ov{b}}^{(1+a+b-n)}\label{eqn:iserre111}
		\\ \notag =&\sum_{n=0,2\mid n}^{a+b+1}\sum_{k=0}^{\frac{n}{2}} \sum_{m=0}^{\lfloor\frac{a+b+1-n}{2}\rfloor} \sqq^{-(1+a+b-n-2m)b} \sum_{r=0}^{\min\{n-2k,1+a+b-n-2m\}}\sum\limits_{[M]\in\cI_{1+a+b-2k-2m-2r}}
		\\ \notag &\quad\quad\frac{\sqq^{\tilde z}(\sqq-\sqq^{-1})^{2+a+b-k-m-r}}{[2k]_\sqq^{!!}[2m]_\sqq^{!!}[r]_\sqq^{!}}\begin{bmatrix} u_M-w_M \\ 1+a+b-n-2m-r-w_M \end{bmatrix}_\sqq
		\frac{[M]*[\E_1]^{k+m+r}}{|\Aut(M)|}
		\\ \notag &-\sum_{n=0,2\nmid n}^{a+b+1}\sum_{k=0}^{\frac{n-1}{2}} \sum_{m=0}^{\lfloor\frac{a+b+1-n}{2}\rfloor} \sqq^{-(1+a+b-n-2m)b} \sum_{r=0}^{\min\{n-2k,1+a+b-n-2m\}}\sum\limits_{[M]\in\cI_{1+a+b-2k-2m-2r}}
		\\ \notag &\quad\frac{\sqq^{\tilde z+2k-2m}(\sqq-\sqq^{-1})^{2+a+b-k-m-r}}{[2k]_\sqq^{!!}[2m]_\sqq^{!!}[r]_\sqq^{!}}\begin{bmatrix} u_M-w_M \\ 1+a+b-n-2m-r-w_M \end{bmatrix}_\sqq \frac{[M]*[\E_1]^{k+m+r}}{|\Aut(M)|}.
	\end{align}
	
	Set
	\begin{align*}
		d=k+m+r.
	\end{align*}
	Now we only need to prove that the coefficient of $\frac{[M]*[\E_1]^d}{|\Aut(M)|}$ in the \text{RHS} of \eqref{eqn:iserre111} is zero, for any given $[M]\in\cI_{1+a+b-2d}$ and any $d\in \mathbb{N}$. For simplicity, we denote by $u=u_M$, $w=w_M$ in the following.

	Note the powers of $(\sqq-\sqq^{-1})$ in all terms are the same (equals to $2+a+b-d$). Denote
	\begin{align}
		\widetilde{T}(a,b,d,u,w)=\sum_{n=0,2\mid n}^{a+b+1}\sum_{k=0}^{\frac{n}{2}} \sum_{m=0}^{\lfloor\frac{a+b+1-n}{2}\rfloor}& \sqq^{-(1+a+b-n-2m)b}\delta\{0\leq r\leq n-2k\}
		\label{eqn:identityT}
		\\ \notag &\cdot\frac{\sqq^{\tilde z}}{[2k]_\sqq^{!!}[2m]_\sqq^{!!}[r]_\sqq^{!!}}\begin{bmatrix} u-w \\ 1+a+b-n-2m-r-w \end{bmatrix}_\sqq
		\\ \notag -\sum_{n=0,2\nmid n}^{a+b+1}\sum_{k=0}^{\frac{n-1}{2}} \sum_{m=0}^{\lfloor\frac{a+b+1-n}{2}\rfloor}& \sqq^{-(1+a+b-n-2m)b}\delta\{0\leq r\leq n-2k\}
		\\ \notag &\cdot\frac{\sqq^{\tilde z+2k-2m}}{[2k]_\sqq^{!!}[2m]_\sqq^{!!}[r]_\sqq^{!!}}\begin{bmatrix} u-w \\ 1+a+b-n-2m-r-w \end{bmatrix}_\sqq,
	\end{align}
	where we set $\delta\{X\}=1$ if the statement $X$ holds and $\delta\{X\}=0$ if $X$ is false. (Note that $r$ is interpreted to be $d-k-m$ in \eqref{eqn:identityT} and \S\ref{subsec:proofidentity} below.) Then the coefficient of $\frac{[M]*[\E_1]^d}{|\Aut(M)|}$ in the \text{RHS} of \eqref{eqn:iserre111} is equal to $(\sqq-\sqq^{-1})^{2+a+b-d}\widetilde{T}(a,b,d,u,w)$. So we have the following.
	\begin{lemma}
		\label{lem:equivalent}
		The identity \eqref{eqn:iserre10} is equivalent to the identity $\widetilde{T}(a,b,d,u,w)=0$, and for any nonnegative integers $a,b,d,u,w$ subject to the constraints
		\begin{align}
			0\leq d\leq (a+b+1)/2,\quad0\leq w\leq b,\quad b+1-2d\leq u\leq 1+a+b-2d,
			\label{eqn:identityTT}
		\end{align}
		where $u,d$ are not both zero. 
	\end{lemma}
	
	In \S\ref{subsec:proofidentity} below, we shall prove the following result.
	\begin{lemma}
		\label{lem:identity=0}
		For any  nonnegative integers $a,b,d,u,w$ satisfying the constraint \eqref{eqn:identityTT}, we have the following identity
		\begin{align}
			\widetilde{T}(a,b,d,u,w)=0.
		\end{align}
	\end{lemma}
	Thanks to Lemma \ref{lem:identity=0}, the formula \eqref{eqn:iserre10} follows.

	The proof of \eqref{eqn:iserre11} is essentially the same as the proof of \eqref{eqn:iserre10}, since the differences on $[S]_{\overline 0}^{(n)}$ versus $[S]_{\overline 1}^{(n)}$ merely lie in the powers of $\sqq$. Going through the same computations, we can get $\widetilde{T_1}$ by changing the power of $\sqq$ in the first summand from $\tilde z$ to $\tilde{z}+2k-2m$ and the power of $\sqq$ in the second summand from $\tilde{z}+2k-2m$ to $\tilde{z}$ in \eqref{eqn:identityT}. The identity $\widetilde{T_1}(a,b,d,u,w)=0$ follows from the identity $\widetilde{T}(a,b,d,u,w)=0$ by using the same argument as in \cite[\S8.3]{LW20a}, which is omitted here.
	
	\subsection{Proof of Lemma \ref{lem:identity=0}}
	\label{subsec:proofidentity}
	The proof given here is modified from the one of \cite[Proposition 8.1]{LW20a}. 
	It is crucial for our purpose to introduce a new variable
	\begin{align*}
		\theta=n+m-k-d
	\end{align*}
	in place of $n$ in \eqref{eqn:identityT}. Hence we have
	\begin{align*}
		\begin{bmatrix} u-w \\ 1+a+b-n-2m-r-w \end{bmatrix}_\sqq =\begin{bmatrix} u-w \\ 1+a+b-2d-\theta-w \end{bmatrix}_\sqq
	\end{align*}
	and
	\begin{align}
		n=\theta-m+k+d\equiv\theta+r \qquad (\mod 2).
		\label{eqn:coresidual}
	\end{align}
	The condition $r\leq n-2k$ in $\widetilde{T}(a,b,d,u,w)$ in \eqref{eqn:identityT} is transformed into the condition $\theta\geq 0$.
	
	By a direct computation we can rewrite the power of $\sqq$ in \eqref{eqn:identityT} as
	\begin{align}
		\overline{z}=&-(1+a+b-n-2m)b+\tilde{z}
		\label{eqn:huajian}
		=\binom{r+1}{2}-2(k-1)m+\tilde L,
	\end{align}
	where
	\begin{align*}
		\tilde L=d(d-1)-a\theta+(u+w+\theta-b)(1+a+b-2d-\theta)-uw+\theta^{2}+1.
	\end{align*}
	(We do not need the precise formula for $\tilde L$ except noting that $\tilde L$ is independent of $k,m,r$, and only depends on $a,b,d,\theta,u,w$.)
	
	By \cite[Lemma 8.3]{LW20a}, for $d\geq1$, the following identity holds
	\begin{align}
		\sum_{\substack{k,m,r\in\N \\ k+m+r=d}}(-1)^r \frac{v^{\tbinom{r+1}{2}-2(k-1)m}}{[r]^{!}[2k]^{!!}[2m]^{!!}}=0.
		\label{eqn:lemma8.3}
	\end{align}
	Hence, for fixed $a,b,d,\theta,u,w$, using \eqref{eqn:coresidual}--\eqref{eqn:lemma8.3}, we calculate that the contribution to the coefficient of $\begin{bmatrix} u-w \\ 1+a+b-2d-\theta-w \end{bmatrix}_\sqq$ in $\widetilde{T}(a,b,d,u,w)$ in \eqref{eqn:identityT}, for $d>0$, is equal to
	\begin{align*}
		&(-1)^\theta \sqq^{\tilde L} \left(\sum_{\substack{k,m,r\in\N,2\mid r\\ k+m+r=d}}(-1)^r \frac{\sqq^{\tbinom{r+1}{2}-2(k-1)m}}{[r]^{!}_\sqq[2k]_\sqq^{!!}[2m]_\sqq^{!!}}+
		\sum_{\substack{k,m,r\in\N,2\nmid r\\ k+m+r=d}}(-1)^r \frac{\sqq^{\tbinom{r+1}{2}-2k(m-1)}}{[r]_\sqq^{!!}[2k]_\sqq^{!!}[2m]_\sqq^{!!}}\right)
		\\&\overset{(*)}{=}(-1)^\theta \sqq^{\tilde L}\sum_{\substack{k,m,r\in\N\\ k+m+r=d}}(-1)^r \frac{\sqq^{\tbinom{r+1}{2}-2(k-1)m}}{[r]_\sqq^{!}[2k]_\sqq^{!!}[2m]_\sqq^{!!}}=0.
	\end{align*}
	Note that the identity $(*)$ above is obtained by switching notation $k\leftrightarrow m$ in the second summand on the \text{LHS} of $(*)$. Therefore, we have obtained that $\widetilde{T}(a,b,d,u,w)=0$, for $d>0$.
	
	It remains to determine the contributions of the terms with $d=0$ to $\widetilde{T}(a,b,0,u,w)=0$ in \eqref{eqn:identityT}, for fixed $a,b,u,w$. Recall from \eqref{eqn:identityTT} $u>0$ when $d=0$. In this case, we have $k=m=r=0$, and a direct computation shows that the power $\bar z$ (see \eqref{eqn:huajian}) can be simplified to be $\bar z=(u+w-b)(1+a+b)-uw+1+(1+2b-u-w)\theta$. Then for $0<u\leq1+a+b$, we have
	\begin{align*}
		\widetilde{T}(a,b,0,u,w)=\sqq^{(u+w-b)(1+a+b)-uw+1}\sum\limits_{\theta\geq0}(-1)^{\theta}\sqq^{(1+2b-u-w)\theta}
		\begin{bmatrix} u-w \\ 1+a+b-\theta-w \end{bmatrix}_\sqq.
	\end{align*}
	By changing the variable $1+a+b-\theta-w$ to $x$, we get
	\begin{align*}
		&\widetilde{T}(a,b,0,u,w)
		\\
		=&(-1)^{1+a+b-w}\sqq^{(1+b)(1+a+b)-w(1+2b-u-w)-uw+1}\sum\limits_{x\geq0}(-1)^{x}\sqq^{(u+w-1-2b)x}
		\begin{bmatrix} u-w \\ x \end{bmatrix}_\sqq.
	\end{align*}
	Using \eqref{eqn:identityTT} and $(u+w-2b-1)\equiv  (u-w-1) \pmod 2$, one can obtain $\widetilde{T}(a,b,0,u,w)=0$ by \cite[Lemma 5.4 (1)]{LW20a}. The proof is completed.

	\section{$\imath$Quantum groups and $\imath$Hall algebras}
	\label{subsec:iQGiH}
	
	Recall $\I_\btau$ from \eqref{eq:ci}.
	Let $\bvs=(\vs_i)\in   (\Q(\sqq)^\times)^{\I}$ be such that $\vs_i=\vs_{\btau i}$ for each $i\in \I$  which satisfies $c_{i, \btau i}=0$. The \emph{reduced Hall algebra associated to $(Q,\btau)$} (or {\em reduced $\imath$Hall algebra}), denoted by $\rMH$, is defined to be the quotient $\Q(\sqq)$-algebra of $\tMHk$ by the ideal generated by the central elements
	\begin{align}
		\label{eqn: reduce1}
		[\E_i] +q \vs_i \; (\forall i\in \I \text{ with } i=\btau i), \text{ and }\; [\E_i]*[\E_{\btau i}] -\sqq^{c_{i,\btau i}}\vs_i\vs_{\btau i}\; (\forall i\in \I \text{ with }i\neq \btau i).
	\end{align}

	\begin{theorem}
		\label{thm:Ui=iHall}
		Let $(Q, \btau)$ be an arbitrary $\imath$quiver. Then there exists a $\Q({\sqq})$-algebra embedding
		\begin{align}
			\label{eqn:psi morphism}
			\widetilde{\Psi}: \tUi_{|v={\sqq}} &\longrightarrow \tMHk,
		\end{align}
		which sends
		\begin{align}
			B_i \mapsto \frac{-1}{q-1}[S_{i}],\text{ if } i\in\ci,
			&\qquad\qquad 	B_{i} \mapsto \frac{{\sqq}}{q-1}[S_{i}],\text{ if }i\notin \ci
			;
			\label{eq:split}
			\\
			\tk_j \mapsto - q^{-1}[\E_j], \text{ if }\btau j=j,
			&\qquad\qquad
			\tk_j \mapsto \sqq^{\frac{-c_{j,\btau j}}{2}}[\E_j],\quad \text{ if }\btau j\neq j.
			\label{eq:extra}
		\end{align}
		Moreover, it induces an embedding $\Psi: \Ui_{|v={\sqq}}{\rightarrow} \rMH$, which sends $B_i$ as in \eqref{eq:split} and  $k_j \mapsto  \vs_{\btau j}^{-1}\sqq^{\frac{-c_{j,\btau j}}{2}}[\E_j], \text{ for } j \in \I\backslash\ci$.
	\end{theorem}

	\begin{proof}
		Recall $[S_i]_{\ov{p}}^{(m)}$, for $i =\tau i$, defined in \eqref{eq:idividedHallodd}--\eqref{eq:idividedHallev} and $[S_i]^{(m)}$, for $i \neq \tau i$, defined in \eqref{eq:dividedpowerHall}. 
		For any $m\in\N$, the map $\widetilde{\Psi}$ in Theorem~\ref{thm:Ui=iHall} sends the $\imath$divided powers $B_{i,\ov p}^{(m)}$ in \eqref{eq:iDPodd}--\eqref{eq:iDPev} for $i =\tau i$ (cf. \cite[Lemma 6.3]{LW20a}) and the divided powers $B_{i}^{(m)}$ in \eqref{eq:DP} for $i \neq \tau i$ to
		\begin{align}
			\widetilde{\Psi}(B_{i,\ov p}^{(m)})
			&= \frac{[S_i]_{\ov p}^{(m)}}{(1-\sqq^2)^m} \;\; \text{ for } \ov{p} \in \Z_2, \quad
			\quad \text{ if } i =\tau i,
			\label{eq:iDPpsi} \\
			\widetilde{\Psi}(B_{i}^{(m)})
			& =
			\begin{cases}
				\frac{[S_i]^{(m)}}{(1-\sqq^2)^m}     & \text{ for } i \in \I_\btau
				\\
				\frac{ \sqq^m [S_i]^{(m)}}{(\sqq^2-1)^m} & \text{ for } i \not \in \I_\btau,
			\end{cases}
			\qquad \text{ if } i \neq \tau i.
			\label{eq:DPpsi}
		\end{align}
		
		To show that $\widetilde{\Psi}$ is a homomorphism, we verify that $\widetilde{\Psi}$ preserves the defining relations \eqref{relation1}--\eqref{relation6} for $\tUi$.
		According to \cite[Lemma 3.9]{LW20a}, the verification of the relations is local and hence is reduced to the rank 1 and rank 2 $\imath$quivers.
		More precisely, the relation \eqref{relation1} follows from \cite[Proposition 9.1]{LW20a}.
		The relation \eqref{relation2} is obvious.
		The relation \eqref{relation3} follows from \cite[Proposition 9.3]{LW20a} by noting that the same result of \cite[Lemma 9.2]{LW20a} also holds for arbitrary quivers; see \cite[Theorem 3.16]{Sch12}.
		The relation \eqref{relation5} follows from \cite[Proposition 5.7]{LW20a}.
		Finally, the relation \eqref{relation6} follows from Theorem \ref{thm: iserre relation} and the same proof of \cite[Proposition 9.4]{LW20a}.
		
		The injectivity of $\widetilde{\Psi}$ follows by using the same proof of \cite[Theorem 9.6]{LW20a}. The proof is completed.
	\end{proof}

	Let $\tCMH$ be the $\Q(\sqq)$-subalgebra (called the {\em composition algebra}) of $\tMHk$ generated by $[S_i]$ and $[\E_i]^{\pm 1}$, for $i\in\I$. Using the same proof of \cite[Corollary 9.9]{LW20a}, we have the following corollary.
	
	\begin{corollary} 
		\label{cor: composition}
		Let $(Q,\btau)$ be an arbitrary $\imath$quiver. Then there exists an algebra isomorphism:
		$
		\widetilde{\Psi}: \tUi_{|v={\sqq}} \stackrel{\cong}{\longrightarrow} \tCMH
		$ 
		given by \eqref{eq:split}--\eqref{eq:extra}.
	\end{corollary}
	
	Recall an oriented cycle of $Q$ is called \emph{minimal} if
	it does not contain any proper oriented cycle. A minimal cycle of length $m$ is called an \emph{$m$-cycle}. 
	An $\imath$quiver $(Q,\btau)$ is called {\em virtually acyclic} if its only possible cycles are $2$-cycles formed by arrows between $i$ and $\tau i$ for $\tau i \neq i \in Q_0$.

	\begin{remark}
		Theorem \ref{thm:Ui=iHall} and Corollary \ref{cor: composition} generalize \cite[Theorem 7.7, Corollary 9.9]{LW20a} by dropping the assumption of $\imath$quivers $(Q,\btau)$ being \emph{virtually acyclic}.

		For the braid group symmetries of $\tUi$ and their $\imath$Hall algebra realization obtained in \cite{LW21}, by using Theorem \ref{thm:Ui=iHall} and Corollary \ref{cor: composition}, one can also drop the  assumption of $\imath$quivers $(Q,\btau)$ being  \emph{virtually acyclic}.
	\end{remark}

	\begin{remark}
		Let $C_n$ ($n\geq2$) be the oriented cyclic quiver with $n$ vertices. With the help of Theorem \ref{thm:Ui=iHall}, we can use the $\imath$Hall algebra of $C_n$ with trivial involution to realize the  $\imath$quantum groups of $\widehat{\mathfrak{sl}}_n$ and $\widehat{\mathfrak{gl}}_n$, especially for the case $n=2$ where $(C_2,\Id)$ is not a virtually acyclic $\imath$quiver; see \cite[Section 10]{LR21}.	
	\end{remark}



\begin{thebibliography}{99}
		
		\bibitem{BW18a} H. Bao and W. Wang,
		{\em Canonical bases arising from quantum symmetric pairs}, Inventiones Math. {\bf 213} (2018), 1099--1177.
		
		\bibitem{BeW18} C. Berman and W. Wang, {\em Formulae of $\imath$-divided powers in ${\mathbf U}_q(\mathfrak{sl}_2)$}, J. Pure Appl. Algebra {\bf 222} (2018), 2667--2702.
		
		\bibitem{Br} T. Bridgeland,
		{\em Quantum groups via Hall algebras of complexes}, Ann. Math. {\bf 177} (2013), 739--759.
		
		\bibitem{CLW18} X. Chen, M. Lu and W. Wang,
		{\em A Serre presentation for the $\imath$quantum gruops}, Transform. Groups {\bf 26} (2021), 827--857, 
		\href{https://arxiv.org/abs/1810.12475}
		{arXiv:1810.12475}
		
		\bibitem{CLW21} X. Chen, M. Lu and W. Wang,
		{\em Serre-Lusztig relations for $\imath$quantum groups}, Commun. Math. Phys. {\bf 382} (2021), 1015--1059,
		\href{http://arxiv.org/abs/2001.03818}
		{arxiv:2001.03818}
		
		\bibitem{Gor2} M. Gorsky,
		{\em Semi-derived and derived Hall algebras for stable categories}, IMRN, Vol. {\bf 2018}, No . 1,  138--159, \href{https://arxiv.org/abs/1409.6798}
		{arXiv:1409.6798}
		
		\bibitem{Gr95} J.A. Green,
		{\em Hall algebras, hereditary algebras and quantum groups}, Invent. Math. {\bf 120} (1995), 361--377.
		
		
		\bibitem{Ko14} S. Kolb,
		{\em Quantum symmetric Kac-Moody pairs}, Adv. Math. {\bf 267} (2014), 395--469.
		
		\bibitem{Let99} G. Letzter,
		{\em Symmetric pairs for quantized enveloping algebras}, J. Algebra {\bf 220} (1999), 729--767.
		
		\bibitem{Let02}
		G. Letzter,
		{\em Coideal subalgebras and quantum symmetric pairs},
		New directions in Hopf algebras (Cambridge), MSRI publications, {\bf 43}, Cambridge Univ. Press, 2002, pp. 117--166.
		
		
		\bibitem{Lu19} M. Lu,
		{\em Appendix A to  \cite{LW19a}},   \href{http://arxiv.org/abs/1901.11446}
		{arXiv:1901.11446}
		
		\bibitem{LP16} M. Lu and L. Peng,
		{\em Semi-derived Ringel-Hall algebras and Drinfeld double}, Adv. Math. {\bf 383} (2021), 107668, \href{https://arxiv.org/abs/1608.03106}
		{arXiv:1608.03106v3}
		
		\bibitem{LR21} M. Lu and S. Ruan,
		{\em $\imath$Hall algebras of weighted projective lines and quantum symmetric pairs}, \href{https://arxiv.org/abs/2110.02575}
		{arXiv:2110.02575}
		
		
		
		\bibitem{LW19a} M. Lu and W. Wang,
		{\em Hall algebras and quantum symmetric pairs I: foundations}, Proc. London Math. Soc. (3) {\bf124} (2022),  1--82. \href{http://arxiv.org/abs/1901.11446}
		{arXiv:1901.11446}
		
		\bibitem{LW20a} M. Lu and W. Wang,
		{\em Hall algebras and quantum symmetric pairs of Kac-Moody type}, Adv. Math. {\bf430} (2023), 109215. \href{http://arxiv.org/abs/2006.06904}
		{arXiv:2006.06904v2}
		
		\bibitem{LW21} M. Lu and W. Wang,
		{\em Braid group symmetries on quasi-split $\imath$quantum groups via $\imath$Hall 
			algeras}, Selecta Math. {\bf28} (2022), Article number 84. \href{https://arxiv.org/abs/2107.06023}
		{arXiv:2107.06023}
		
		
		\bibitem{Lus93} G. Lusztig, Introduction to Quantum Groups, Birkh\"{a}user, Boston, 1993.
		
		\bibitem{Rin90} C.M. Ringel,
		{\em Hall algebras and quantum groups}, Invent. Math. {\bf 101} (1990), 583--591.
		
		\bibitem{Sch12} O. Schiffmann, {\em Lectures on Hall algebras}, In: Geometric methods in representation theory. II, 1--141,
		S\'{e}min. Congr., 24-II, Soc. Math. France, Paris (2012).
		
	\end{thebibliography}
\end{document}